\documentclass[12pt]{article}

\usepackage{amssymb, amsmath, amsthm}   % For math typesetting
\usepackage{graphicx}  % For including images

\usepackage{hyperref}  % For clickable links
\usepackage{geometry}  % Page layout
\usepackage{multirow}
\usepackage[normalem]{ulem}
\usepackage{hyperref}
\hypersetup{
    colorlinks=true,
    linkcolor=blue,
    filecolor=magenta,
    urlcolor=cyan,
    citecolor=red
}

\newcommand{\comment}[1]{}
\newtheorem{theorem}{\textbf{Theorem}}
\newtheorem{lemma}{\textbf{Lemma}}
\newtheorem{proposition}{\textbf{Proposition}}

\newtheorem{remark}{\textbf{Remark}}

\usepackage{placeins} %keep floats into place

\begin{document}
%\begin{linenumbers}
\begin{center}  {\Large \bf ODE and PDE models for COVID-19, with reinfection and vaccination process for Cameroon and Germany}
	\end{center}
	
	\begin{center}
		{%\footnotesize
\textsc{Hamadjam Abboubakar}$^{\dagger,\mp}$,
\textsc{Reinhard Racke}$^{\pm}$,
\textsc{Nicolas Schlosser}$^{*,\pm}$,
 }
\end{center}
\begin{center}
{\small \sl $^{\dagger}$
Department of Computer Engineering, University Institute of Technology of The University of Ngaoundéré, P.O. Box 455 Ngaoundéré, Cameroon\\}
{\small \sl $^{\mp}$
Department of Applied Mathematics and Computer Science, School of Geology and Mining Engineering of The University of Ngaoundéré, P.O. Box 115 Meiganga, Cameroon}\\
{\small \sl $^{\pm}$
Department of Mathematics and Statistics, University of Konstanz, P.O. Box 78457 Konstanz,  Germany
}\\
\end{center}

	 %\maketitle
\begin{abstract}
The goal of this work is to develop and analyze a reaction-diffusion model for the transmission dynamics of the Coronavirus (COVID-19) that accounts for reinfection and vaccination, as well as to compare it to the ODE model. After developing a time-dependent ODE model, we calculate the control reproduction number $\mathcal{R}_c$ and demonstrate the global stability of the COVID-19 free equilibrium for $\mathcal{R}_c<1$. We also show that when $\mathcal{R}_c>1$, the free equilibrium of COVID-19 becomes unstable and co-exists with at least one endemic equilibrium point. We then used data from Germany and Cameroon to calibrate our model and estimate some of its characteristics. We find $\mathcal{R}_c\approx 1.13$ for Germany and $\mathcal R_c \approx  1.2554$ for Cameroon, indicating that the disease persists in both populations. Following that, we modify the prior model into a reaction-diffusion PDE model to account for spatial mobility. We show that the solutions to the final initial value boundary problem (IVBP) exist and are nonnegative and unique. We also show that the disease-free equilibrium is stable locally, and globally when $\mathcal{R}_c<1$. In contrast, when $\mathcal{R}_c>1$, the DFE is unstable and coexists with at least one endemic equilibrium point. We ran multiple numerical simulations to validate our theoretical predictions. We then compare the ODE and the PDE models.


\end{abstract}

\section{Introduction}
\label{sec1}
The Covid-19 pandemic started in 2019, see \cite{world2020coronavirus} for development and containment measurements. May 2021 is the start date for vaccination in several countries around the world \cite{WhoVaccine}.
Since the beginning of the Covid-19 pandemic, several mathematical models have been formulated and studied to predict the future of the disease, as well as the efficiency of control measures (see \cite{adiga2020mathematical,clement2021survey,harjule2021mathematical,james2021use,nabi2020forecasting,rahimi2021review} and the references therein). Some authors have worked on time-space models, also called reaction-diffusion models \cite{ahmed2021numerical,brusset2021modelling,cherniha2021reaction,fitzgibbon2020analysis,grave2022modeling,kevrekidis2021reaction,mammeri2020reaction,turkyilmazoglu2022indoor,zhu2021dynamic}.
In \cite{ahmed2021numerical}, Ahmed et {\it al.} formulated a SAIR reaction-diffusion model with nonlinear incidence rates in a constant population. Brusset et {\it al.} in \cite{brusset2021modelling} formulated a SIS reaction-diffusion model to represent how the geographic spread of the pandemic, by reducing the workers' participation to economic life,  undermines the ability of firms and as a result the entire supply networks to satisfy customers' demands. An SI reaction-diffusion model with cross-diffusion is formulated and studied in \cite{cherniha2021reaction} by Cherniha and Davydovych using the Lie symmetry method. In \cite{fitzgibbon2020analysis}, Fitzgibbon et {\it al.} developed a dynamic model of an evolving epidemic in a spatially inhomogeneous environment. They analyzed it to predict the outbreak and spatio-temporal spread of the COVID-19 epidemic in Brazil. To take into account the nonlocal Covid-19 transmission due to the fact that people often travel long distances in short periods of time, Grave et {\it al.} in \cite{grave2022modeling} combined a network structure within a reaction-diffusion PDE system. They defined the transfer network, the transfer operator, the donor operator, and the receiver operator.
Kevrekidis et {\it al.} \cite{kevrekidis2021reaction} formulated and studied an Susceptible-Exposed-Asymptomatic-Infectious-Hospitalized-Recovered (SEAIHR) reaction-diffusion model with Greece and Andalusia as case examples. Youcef Mammeri in \cite{mammeri2020reaction} formulated and studied a SEAIR reaction-diffusion model with mass action incidences, and France as a case example. Mustafa Turkyilmazoglu in \cite{turkyilmazoglu2022indoor} formulated and studied a simplistic reaction-diffusion model to mathematically explore the spatiotemporal development of the concentration of indoor aerosols containing infectious nuclei of the COVID-19 respiratory virus. Zhu and Zhu in \cite{zhu2021dynamic} constructed a time delay reaction-diffusion model that includes relapse, time delay, home quarantine, and a heterogeneous temporal-spatial environment that affect the spread of COVID-19.

Note that the authors mentioned above do not integrate vaccination in their models. In \cite{Hamadjam-Reinhard}, we formulated and studied a Susceptible-Quarantined-Vaccinated-Exposed-Asymptomatic infectious-symptomatic infectious-Hospitalized-Recovered- Virus\\ ($SQVEAIHR-B$) Covid-19 type compartmental model in which vaccinated individuals are divided into two different groups: those who take the first dose and those who take the second dose after taking the first dose. Model formulation was done using both integer and noninteger derivative in the Caputo Sense, with application to German data. In the present study, we formulate and study a reaction-diffusion model of type Susceptible-Vaccinated-Exposed-Asymptomatic Infectious- Symptomatic Infectious-Recovered (SVEAIR) to translate the transmission dynamics of Covid-19. The model proposed here takes into account the reinfection and vaccination process with Germany  or Cameroon as a case study. One main goal here is to compare quantitatively the model formulated with ordinary differential equations (ODE) and the corresponding reaction-diffusion model with partial differential equation (PDE).
Our new contributions are:
\begin{enumerate}
	\item We first analyze the ODE model by determining the control reproduction number denoted by $\mathcal{R}_c$ and prove the global asymptotic stability of the disease-free equilibrium point (DFEP) whenever $\mathcal{R}_c<1$. Then, we prove the existence of at least one endemic equilibrium point (EEP) when $\mathcal{R}_c>1$.
	\item We then perform parameter estimation using real data from Germany and Cameroon.
	\item After that, we extend the ode model by including the diffusion terms to obtain a reaction-diffusion PDE model. We prove the non-negativity of state variables as well as the existence and uniqueness of solutions. The asymptotic stability results of the DFEP of the ODE model are extended to obtain the asymptotic stability of the DFEP of the PDE model whenever the control reproduction number $\mathcal{R}_c$ is less than one.
	\item Numerical simulations are finally performed by considering:
	\begin{enumerate}
		\item Constant parameters and time-dependent parameters,
		\item Two cases:
		\begin{enumerate}
			\item Initial population is completely susceptible to infection everywhere except for one small region in the very south of Germany, where there are also infected persons;
			\item  When we add a second peak in Germany -- resp. Cameroon --, where a major outbreak of COVID-19 occurred in early 2020 -- resp. in March for Cameroon.
		\end{enumerate}
	\end{enumerate}
\end{enumerate}
In the present study, we have not yet considered that the model parameters can be time-dependent or space-dependent up to the numerical treatment of a time-dependent transmission rate $\beta$.

 The paper revises and extends -- to the study of Cameroon -- the preliminary version \cite{abboubakar2023reaction}.

The outline of the work is as follows: the formulation of the ODE compartmental model and its theoretical analysis is given in Section \ref{Model_formulation}. Section \ref{calibration} is devoted to the calibration, forecasting, and global sensitivity analysis of the model. The reaction-diffusion model is formulated and studied in Section \ref{Reaction-Diffusion}. Section \ref{num_analysis} is devoted to the numerical scheme and the simulation results.

\section{Model formulation}
\label{Model_formulation}
The model we consider here is an extension of an SEIR-type compartmental model, in which we take into account reinfection as well as the vaccination process. The total population at each time $t$, denoted by $\mathfrak{N}(t)$, is divided into six states or compartments as follows: susceptible people denoted by $\mathfrak{S}(t)$, vaccinated people denoted by $\mathfrak{V}(t)$, infected people in the latent stage denoted by $\mathfrak{E}(t)$, infected people without symptoms (asymptomatic) denoted by $\mathfrak{A}(t)$, infected people with symptoms (symptomatic) denoted by $\mathfrak{I}(t)$, and recovered people denoted by $\mathfrak{R}(t)$. So, $\mathfrak{N}(t)=\mathfrak{S}(t)+\mathfrak{V}(t)+\mathfrak{E}(t)+\mathfrak{A}(t)+\mathfrak{I}(t)+\mathfrak{R}(t)$. In this model,  the compartment $\mathfrak{I}(t)$ includes all detected cases as well as hospitalized cases, while $\mathfrak{A}(t)$ includes all people who are infectious, but are not tested and do not have any symptoms of the disease. We consider immigration of vaccinated people into the system. Thus, considering the parameter $\Lambda$ as the recruitment rate of non-infected people, a rate $r_2$ of these recruited people is vaccinated while the rest denoted by $r_1$ is not vaccinated. Among these non-vaccinated people, a rate of $c_1$ will be vaccinated while among vaccinated people, a rate of $c_2$ will lose their immunity conferred by the vaccine and become again susceptible. Susceptible people can contract the virus by direct contact with asymptomatic or symptomatic individuals at the rate $\lambda(t)=\dfrac{\beta\left(\mathfrak{A}(t)+\eta \mathfrak{I}(t)\right)}{\mathfrak{N}(t)}$, where $\beta$ represents the transmission rate, while $\eta$ represents the modification parameter due to the fact that people who have been tested positive are considered less infectious because they must take control measures (isolation, quarantine, treatment, \ldots) to limit the transmission of the disease. To represent the efficacy of the vaccine, we introduce the parameter $\phi_{1}=(1-\epsilon)$, where $\epsilon$ represents the efficacy of the Covid-19 vaccine. Thus, the fraction of $\phi_{1}\lambda(t)(t)$ will become infected after close contact with an infectious individual. After $1/\gamma$ days, which represents the latent period, infected people will become either in the $A$ compartment or in the $I$ compartment. Asymptomatic people will move either in the compartment $I$ or in the recovered compartment $R$ at the rates $a_{2}\sigma$ and $a_{1}\sigma$, respectively. Symptomatic people can recover from infection naturally or after treatment at a rate $\theta$. In the opposite case, some of them will die of the disease. So, the parameter $\delta$ represents the disease-induced death. Recovered people can become reinfected at a rate $\phi_{2}\lambda(t)(t)$, where $\phi_{2}$ represents the rate of recovered people who will become infected again. In each compartment, people can die naturally with a natural death rate $\mu$.

The COVID-19 transmission dynamics model expressed using ODEs looks as follows:
\begin{equation}
	\label{Covid19_model-reduce}
	\left\lbrace \begin{array}{ll}
		\dfrac{d\mathfrak{S}(t)}{dt}&=r_1\Lambda+c_2\mathfrak{V}(t)-\left[\overbrace{c_1+\mu}^{k_1}
		+\overbrace{\dfrac{\beta\left( \mathfrak{A}(t)+\eta \mathfrak{I}(t)\right)}{\mathfrak{N} (t)}}^{\lambda(t)(t)}\right] \mathfrak{S}(t),\\
		\dfrac{d\mathfrak{V}(t)}{dt}&=r_2\Lambda+c_1\mathfrak{S}(t)-\left[\overbrace{\mu+c_2}^{k_2}+\phi_1\dfrac{\beta\left( \mathfrak{A}(t)+\eta \mathfrak{I}(t)\right)}{\mathfrak{N}(t)} \right] \mathfrak{V}(t),\\
		\dfrac{d\mathfrak{E}(t)}{dt}&=\dfrac{\beta\left( \mathfrak{A}(t)+\eta \mathfrak{I}(t)\right)}{\mathfrak{N}(t)}\left(\mathfrak{S}(t)+\phi_1\mathfrak{V}(t)+\phi_2\mathfrak{R}(t)\right) -\overbrace{\left(\mu+\gamma\right)}^{k_3}\mathfrak{E}(t) ,\\
		\dfrac{d\mathfrak{A}(t)}{dt}&=p\gamma \mathfrak{E}(t)-\overbrace{(\mu+\sigma)}^{k_4}\mathfrak{A}(t),\\
		\dfrac{d\mathfrak{I}(t)}{dt}&=q\gamma \mathfrak{E}(t)+a_2\sigma \mathfrak{A}(t)-\overbrace{\left( \mu+\delta+\theta\right) }^{k_5}\mathfrak{I}(t),\\
		\dfrac{d\mathfrak{R}(t)}{dt}&=a_1\sigma \mathfrak{A}(t)+\theta \mathfrak{I}(t)-\left[ \mu+\phi_2\dfrac{\beta\left( \mathfrak{A}(t)+\eta \mathfrak{I}(t)\right)}{\mathfrak{N}(t)}\right]\mathfrak{R}(t)\\
	\end{array} \right.
\end{equation}
with $r_1+r_2=1$, $p+q=1$, and $a_1+a_2=1$.

Setting $X=\left(\mathfrak{S},\mathfrak{V},\mathfrak{E},\mathfrak{A},\mathfrak{I},\mathfrak{R}\right)'$ the vector of state variables\\ and
$\Pi=\left\lbrace
X\in\mathbb{R}^{6}:X\geq \mathbf{0}_{\mathbb{R}^{6}}\right\rbrace$, system \eqref{Covid19_model-reduce} can be written in the following compact form
\begin{equation}
	\label{Covid19_model-reduce-compact}
	\left\{ \begin{array}{l}
		\dfrac{dX}{dt}=\mathcal{F}(t,X)=\left(\mathcal{F}_{1}(X),\mathcal{F}_{2}(X),...,\mathcal{F}_{6}(X)\right)',\\
		X(t_0)=X_{0}=(\mathfrak{S}_0,\mathfrak{V}_{0},\mathfrak{E}_0,\mathfrak{A}_0,\mathfrak{I}_0,\mathfrak{R}_0)'\geq \mathbf{0}_{\mathbb{R}^{6}},
	\end{array}
	\right.
\end{equation}
where $\mathcal{F}:\mathbb{R}^{6}\to\mathbb{R}^{6}$ is a continuously differentiable function on $\mathbb{R}^{6}$, and $(\bullet)'$ stands for the transposition operator. According to \cite[Theorem III.10.VI]{walter2013ordinary}, for $X(0)\in \Pi$, a unique solution of \eqref{Covid19_model-reduce} exists, at least locally, and remains in $\Pi$ for its maximal interval of existence \cite[Theorem III.10.XVI]{walter2013ordinary}. Hence model \eqref{Covid19_model-reduce} is biologically well-defined.

Model \eqref{Covid19_model-reduce} is defined in the following set
\[
\mathbf{W}=\left\lbrace X:=\left(\mathfrak{S},\mathfrak{V},\mathfrak{E},\mathfrak{A},\mathfrak{I},\mathfrak{R}\right)'\in\mathbb{R}^{6}_{+}:0<N:=\sum\limits_{i=1}^{6}X_i\leq\dfrac{\Lambda}{\mu}\right\rbrace,
\]
which is invariant for the system \eqref{Covid19_model-reduce}.	

The above statement is obtained in the same way as the results obtained in \cite[Theorem~2]{Hamadjam-Reinhard}.

\subsection{The disease-free equilibrium and the control reproduction number}
In the absence of disease, i.e. for $A=I=B=0$, model \eqref{Covid19_model-reduce} always admits the equilibrium $\mathcal E_0=\left(\mathfrak{S}_0,\mathfrak{V}_{0},0,0,0,0\right)'$ called the disease-free equilibrium, with \begin{equation}\label{DFE}
	\mathfrak{S}_0 =\frac{\left(c_{2}r_{2}+r_{1}k_{2}\right)\Lambda}{k_{1}k_{2}-c_{1}c_{2}} \quad \text{and}\quad \mathfrak{V}_0 =\frac{\left(k_{1}r_{2}+c_{1}r_{1}\right)\Lambda}{k_{1}k_{2}-c_{1}c_{2}}.
\end{equation} Note that
$k_{1}k_{2}-c_{1}c_{2}=\mu^2+\left(c_{2}+c_{1}\right)\mu>0$, and $\mathfrak{S}_0+\mathfrak{V}_0=\mathfrak{N}_0=\dfrac{\Lambda}{\mu}$.

To compute the control reproduction number, denoted by $\mathcal R_c$, we will use the next generation approach (see \cite{diekmann1990definition,van2002reproduction}).
Let us set $y=\left(E,A,I\right)'$. The vectors $\mathcal{Z}$ and $\mathcal{W}$ for the new infection terms and the remaining transfer terms for $y$ are, respectively, given by\\
$
\mathcal{Z}=\left(\begin{array}{c}
	\dfrac{\beta\left( \mathfrak{A}(t)+\eta \mathfrak{I}(t)\right)}{N(t)}\left(S+\phi_1V+\phi_{2}R\right) \\
	0\\
	0\\
\end{array}\right)
$
and
$
\mathcal{W}=\left(\begin{array}{c}
	k_3E\\
	-p\gamma E+k_4A,\\
	-q\gamma E-a_2\sigma A+k_5I,\\
\end{array}\right)
$.\\
Their Jacobian matrices evaluated at $\mathcal{E}_0$ are respectively given by
\begin{equation}
	\label{ZW_Rc}
	\begin{split}
		\mathfrak{Z}=\left(\begin{array}{ccc}
			0&\beta \dfrac{N_{1}}{N_{0}}&\beta\eta \dfrac{N_{1}}{N_{0}}\\
			0&0&0\\
			0&0&0\\
		\end{array}\right)\,\,\text{and}\,\,\,
		\mathfrak{W}=\left(\begin{array}{ccc}
			k_3&0&0\\
			-p\gamma&k_4&0,\\
			-q\gamma&-a_2\sigma&k_5,\\
		\end{array}\right),
	\end{split}
\end{equation}
with $\mathfrak{N}_{1}=\mathfrak{S}_{0}+\phi_1\mathfrak{V}_{0}$.
Then, the control reproduction number $\mathcal{R}_c$ is defined, following \cite{diekmann1990definition,van2002reproduction}, as the spectral radius of the next generation matrix, $ZW^{-1}$ where
\[
\mathfrak{Z}\mathfrak{W}^{-1}=\left(\begin{array}{ccc}
	\dfrac{N_{1}\beta\eta\left(a_{2}p\sigma
		\gamma+k_{4}q\gamma\right)}{N_{0}k_{3}k_{4}k_{5}}+\frac{N
		_{1}\beta p\gamma}{N_{0}k_{3}k_{4}} & \dfrac{N_{1}a_{2}
		\beta\eta\sigma}{N_{0}k_{4}k_{5}}+\dfrac{N_{1}\beta}{N_{0}
		k_{4}} & \dfrac{N_{1}\beta\eta}{N_{0}k_{5}} \\ 0 & 0 & 0 \\ 0
	& 0 & 0 \\
\end{array}\right).
\]
Therefore, the control reproduction number, $\mathcal{R}_{c}$, is the sum of two main contributions, namely, humans and environment, as follows:
\begin{equation}
	\begin{split}
		\label{Rc}
\mathcal{R}_{c}&:=\rho(\mathfrak{Z}\mathfrak{W}^{-1})= \dfrac{\mathfrak{N}_{1}\beta\eta\gamma\left(a_{2}p\sigma+k_{4}q\right)}{\mathfrak{N}_{0}k_{3}k_{4}k_{5}}+\dfrac{\mathfrak{N}_{1}\beta p\gamma}{\mathfrak{N}_{0}k_{3}k_{4}},
	\end{split}
\end{equation}
where $\rho(\bullet)$ represents the spectral radius operator.

From \cite[Theorem 2]{van2002reproduction}, we have the following result.
\begin{lemma}(Local stability of the DFE)
	\label{LAS_dfe}
	The stationary point $\mathcal{E}_0$ of system \eqref{Covid19_model-reduce} is locally asymptotically stable (LAS) if $\mathcal{R}_{c}<1$, and unstable otherwise.
\end{lemma}
We also have the following result:
\begin{theorem}
	\label{gas_dfe_reduceModel}
	The disease-free equilibrium $\mathcal{E}_0$ is globally asymptotically stable in $\mathcal W$ whenever $\mathcal{R}_{c}<1$, provided that
    \begin{equation}
        \label{cond_gas_dfe}
        \dfrac{N_{1}}{N_{0}}-\dfrac{\left(\mathfrak{S}(t)+\phi_1\mathfrak{V}(t)+\phi_2\mathfrak{R}(t)\right)}{N(t)}\geq 0.
    \end{equation}
\end{theorem}
\begin{proof}	
	Considering only the infected compartments of system \eqref{Covid19_model-reduce}, we obtain
	\begin{equation}
		\left(\begin{array}{c}
			\dfrac{d\mathfrak{E}(t)}{dt}\\\dfrac{d\mathfrak{A}(t)}{dt}\\\dfrac{d\mathfrak{I}(t)}{dt}
		\end{array}\right)=\left(Z-W\right)\left(\begin{array}{c}\mathfrak{E}(t)\\\mathfrak{A}(t)\\\mathfrak{I}(t)\end{array}\right)
		-\mathcal{M}\left(\mathfrak{S}(t),\mathfrak{V}(t),\mathfrak{E}(t),\mathfrak{A}(t),\mathfrak{I}(t),\mathfrak{R}(t)\right),
	\end{equation}
	where $Z$ and $W$ are the same matrices used to compute the control reproduction number (see Eq. \eqref{Rc}), and
	\[
    \begin{split}
&\mathcal{M}(\mathfrak{S}(t),\mathfrak{V}(t),\mathfrak{E}(t),\mathfrak{A}(t),\mathfrak{I}(t),\mathfrak{R}(t))\\
&=\begin{pmatrix}
		\beta\left(\mathfrak{A}(t)+\eta \mathfrak{I}(t)\right) \left(\dfrac{N_{1}}{N_{0}}-\dfrac{\left(\mathfrak{S}(t)+\phi_1\mathfrak{V}(t)+\phi_2\mathfrak{R}(t)\right)}{N(t)} \right)\\
		0\\0 \end{pmatrix}.
    \end{split}
	\]
	%In $\mathcal W$, $N_{0}\geq N_1:=(S_0+\phi_1V_0)\geq\left(S+\phi_1V+\phi_2R\right)$ for all $t>0$.
	If $\dfrac{N_{1}}{N_{0}}-\dfrac{\left(\mathfrak{S}(t)+\phi_1\mathfrak{V}(t)+\phi_2\mathfrak{R}(t)\right)}{N(t)}\geq 0$ in $\mathcal W$,	
	then, it follows that\\ $\mathcal{M}(S(t),V(t),E(t),A(t),I(t),R(t))\geq \mathbf{0}_{\mathbb{R}^{3}}$. This means that
	\[
	\left(\begin{array}{c}
		\dfrac{dE(t)}{dt}\\\dfrac{dA(t)}{dt}\\\dfrac{dI(t)}{dt}
	\end{array}\right)\leq\left(Z-W\right)\left(\begin{array}{c}\mathfrak{E}(t)\\\mathfrak{A}(t)\\\mathfrak{I}(t)\end{array}\right).
	\]
	Note that $W^{-1}=\begin{pmatrix}\frac{1}{k_{3}} & 0 & 0 \\ \frac{p\gamma}{k_{3}k_{4}} & \frac{1}{k_{4}} & 0 \\ \frac{a_{2}p\sigma\gamma+k_{4}
			q\gamma}{k_{3}k_{4}k_{5}} & \frac{a_{2}\sigma}{k_{4}k_{5
		}} & \frac{1}{k_{5}} \\ \end{pmatrix}\geq \mathbf{0}_{\mathbb{R}^{3\times3}}$.\\
	We also have, from \eqref{ZW_Rc}, $Z\geq 0$. Thus, from \cite[Theorem 2.1]{shuai2013global}, there exists a Lyapunov function for system \eqref{Covid19_model-reduce} expressed as $\mathcal{Q}\left(S,V,E,A,I,R\right) =w'W^{-1}\left(E,A,I\right)'$ where $w'$ is the left eigenvector of the nonnegative matrix $W^{-1}Z$ corresponding to the eigenvalue $\mathcal{R}_c$. This implies that if $\mathcal{R}_c<1$, \[
	\dfrac{d\mathcal{Q}}{dt}=\left(\mathcal{R}_c-1\right)w'\left(E,A,I\right)-w'W^{-1}\mathcal{M}\left(S,V,E,A,I,R\right)\leq 0
	\] whenever the condition  $\dfrac{N_{1}}{N_{0}}-\dfrac{\left(\mathfrak{S}(t)+\phi_1\mathfrak{V}(t)+\phi_2\mathfrak{R}(t)\right)}{N(t)}\geq 0$ holds. It follows from the LaSalle invariance principle \cite{la1976stability} that every solution of \eqref{Covid19_model-reduce} with initial conditions in $\mathbf{W}$ converges to the DFE $\mathcal{E}_0$ when $t\longrightarrow+\infty$. That is $\left(E,A,I\right)\longrightarrow\left(0,0,0\right)$, $S\longrightarrow S_{0}$ and $V\longrightarrow V_{0}$ when $t\longrightarrow+\infty$, which is equivalent to $\left(S,V,E,A,I,R\right)\longrightarrow\left(S_{0},V_0,0,0,0,0\right)$ when $t\longrightarrow+\infty$. Thus, by the LaSalle invariance principle \cite{la1976stability}, the disease-free equilibrium $\mathcal E_0$ is globally asymptotically stable in $\mathcal{W}$ whenever $\mathcal{R}_{c}<1$. This ends the proof of Theorem \ref{gas_dfe_reduceModel}.
\end{proof}
%\begin{remark}
%Note that in the case of perfect vaccine, that is $\phi_1=0$, condition \eqref{cond_gas_dfe} always holds since $\dfrac{N_{1}}{N_{0}}-\dfrac{\left(\mathfrak{S}(t)+\phi_2\mathfrak{R}(t)\right)}{N(t)}\geq 1-\dfrac{\left(S(x,t)+\phi_2R(x,t)\right)}{N(t)}$
%\end{remark}

\subsection{Existence of the endemic equilibrium}
Let $\mathcal E=\left(S^{*},V^{*},E^{*},A^{*},I^{*},R^{*}\right)'$ be an equilibrium point of model \eqref{Covid19_model-reduce} obtained by setting the right hand-side of \eqref{Covid19_model-reduce} equal to zero, that is
\begin{equation}
	\label{ee1}
	\left\lbrace \begin{array}{ll}
		r_1\Lambda+c_2V^{*}-k_1S^{*}-\lambda^{*}S^{*}&=0,\\
		r_2\Lambda+c_1S^{*}-\left[k_2+\phi_1\lambda^{*}\right] V^{*}&=0,\\
		\lambda^{*}\left(S^{*}+\phi_1V^{*}+\phi_2R^{*}\right)-k_3E^{*}&=0 ,\\
		p\gamma E^{*}-k_4A^{*}&=0,\\
		q\gamma E^{*}+a_2\sigma A^{*}-k_5I^{*}&=0,\\
		a_1\sigma A^{*}+\theta I^{*}-\left[ \mu+\phi_2\lambda^{*}\right]R^{*}&=0.\\
	\end{array} \right.
\end{equation}
Solving the above system gives
\begin{equation}
	\label{ee2}
	\left\lbrace \begin{array}{ll}
		S^{*}&=\dfrac{r_{1}\Lambda \phi_{1}\lambda^{*}+\overbrace{\left(c_{2}r_{2}+r_{1}k_{2}\right)}^{k_7}\Lambda}
		{\lambda^{*}\left( \phi_{1}\lambda^{*}+k_{2}\right)
			+k_{1}\phi_{1}\lambda^{*}+\underbrace{k_{1}k_{2}-c_{1}c_{2}}_{k_6}}
		=\dfrac{r_{1}\Lambda \phi_{1}\lambda^{*}+k_7\Lambda}
		{\lambda^{*}\left( \phi_{1}\lambda^{*}+k_{2}\right)+k_{1}\phi_{1}\lambda^{*}+k_6},\\
		V^{*}&=\dfrac{r_2\Lambda+c_1S^{*}}{\left[k_2+\phi_1\lambda^{*} \right]},
		E^{*}=\dfrac{\lambda^{*}\left(S^{*}+\phi_1V^{*}+\phi_2R^{*}\right)}{k_3},
		A^{*}=\dfrac{p\gamma E^{*}}{k_4},\\
		I^{*}&=\dfrac{q\gamma E^{*}+a_2\sigma A^{*}}{k_5},
		R^{*}=\dfrac{a_1\sigma A^{*}+\theta I^{*}}{\left[ \mu+\phi_2\lambda^{*}\right]},
		N^{*}=\dfrac{\Lambda-\delta I^{*}}{\mu},
	\end{array} \right.
\end{equation}
where $\lambda^{*}:=\beta\dfrac{\left(A^{*}+\eta I^{*}\right)}{N^{*}}$
is any nonnegative solutions of the following equation
\begin{equation}
	\label{Pol_ee}
	\lambda^{*}\left[\mathcal A_{3}(\lambda^{*})^3+\mathcal A_{2}(\lambda^{*})^2+\mathcal A_{1}\lambda^{*}+\mathcal A_{0} \right]=0,
\end{equation}
where
\begin{align*}
	\mathcal A_{3}&=-\phi_{1}\phi_{2}\left(r_{1}\mu\left(1-\phi_{1}\right)+\phi_{1}\left( \mu+c_{1}\right) +
	c_{2}\right)\left(a_{2}\eta
	p\sigma+k_{4}\eta q+k_{5}p\right)\\
	&\quad \times\left[ \left(\theta-\delta \right) \left( a_{2}p\sigma\gamma+k_{4}q\gamma\right)
	+k_{5}\mu(\mu+k_{2})+k_{5}\gamma\sigma(1-a_1p)\right],\\
	\mathcal A_{0}&=k_{3}k_{4}k_{5}\mu^2\left(\mu+c_{2}+c_{1}\right)\left(r_{1}\mu\left(1-\phi_{1}\right)+\phi_{1}\left( \mu+c_{1}\right)+c_{2}\right)\left(a_{2}\eta p\sigma+k_{4}\eta q+k_{5}p\right)\\
	&\times\left(\mathcal{R}_{c}-1\right),
\end{align*}
Note that the signs of the coefficients $\mathcal A_{2}$ and $\mathcal A_{1}$ depend on the values of the parameters, and that their expressions can be obtained with any formal calculation software such as Maxima$^\copyright$ and Maple$^\copyright$. As these expressions are very long, we simply omit their expression.
Assume that $\theta_1>\delta$. Since $a_1+a_2=1$ and $p+q=1$, it follows that $1-a_1p>0$. Thus, the coefficient $\mathcal A_{3}$ is always negative. Coefficient $\mathcal A_{0}$ is negative (resp. positive) if and only if $\mathcal R_c<1$ (resp. $\mathcal R_c>1$). Using the Descartes rule of signs, we claim the following:
\begin{proposition}$\;$
\label{existence_EE}
\begin{enumerate}
	\item If $\mathcal R_c>1$, then model \eqref{Covid19_model-reduce} admits:
	\begin{itemize}
		\item[(i)] exactly one endemic equilibrium point if and only if
		$\left(\mathcal A_{2}>0\,\,\&\,\,\mathcal A_{1}>0\right)$ or\\ $\left(\mathcal A_{2}<0\,\,\&\,\,\mathcal A_{1}>0\right)$ or $\left(\mathcal A_{2}<0\,\,\&\,\,\mathcal A_{1}<0\right)$;
		\item[(ii)] exactly three endemic equilibrium points iff
		$\left(\mathcal A_{2}>0 \,\&\,\mathcal A_{1}<0\right)$;
	\end{itemize}
	\item If $\mathcal R_c<1$, then model \eqref{Covid19_model-reduce} admits:
	\begin{itemize}
		\item[(iii)] exactly two endemic equilibrium points if and only if
		$\left(\mathcal A_{2}>0\,\,\&\,\,\mathcal A_{1}>0\right)$ or\\ $\left(\mathcal A_{2}>0\,\,\&\,\,\mathcal A_{1}<0\right)$ or $\left(\mathcal A_{2}<0\,\,\&\,\,\mathcal A_{1}>0\right)$;
		\item[(iv)] Otherwise no endemic equilibrium exists.
	\end{itemize}
\end{enumerate}
\end{proposition}
The item $(iii)$ of the above proposition suggests the possibility of the occurrence of the backward bifurcation phenomenon \cite{abboubakar2018bifurcation} in the model \eqref{Covid19_model-reduce}, i.e., when the disease-free equilibrium point coexists with two endemic equilibrium points (one is locally stable and the other is unstable) whenever the biological threshold $\mathcal{R}_c$ is less than one. Thanks to Theorem \ref{gas_dfe_reduceModel}, we conclude that even if the disease-free equilibrium coexists with two endemic equilibrium points when $\mathcal{R}_c<1$, these last ones are either always unstable or do not belong to the set $\mathcal W$.

\section{Model calibration and forecasting}
\label{calibration}
\subsection{Model calibration with German data}
The start date for mass vaccination in Germany was Sunday 27 December 2020 \cite{franceInfo}. Since then, several constraint measures have been taken to ensure that the majority of inhabitants are vaccinated. We consider the daily reported cases of infection in Germany from December 31, 2020, to February 28, 2021 \cite{DataWorldCovid}.

In model \eqref{Covid19_model-reduce}, there are 18 parameters. Of these, seven parameters are either estimated or taken from the existing literature, while the other remaining parameters must be calibrated using data. Taking the approximate total population of Germany in 2021 equal to $N(0)=83,900,473$ \cite{ganegoda2021interrelationship,sitePopulGerman}, the recruitment rate is equal to $\Lambda=\mu N(0)$. %At the December 27, 20201,

\begin{figure}[ht]
\begin{center}
	\includegraphics[width=1\textwidth]{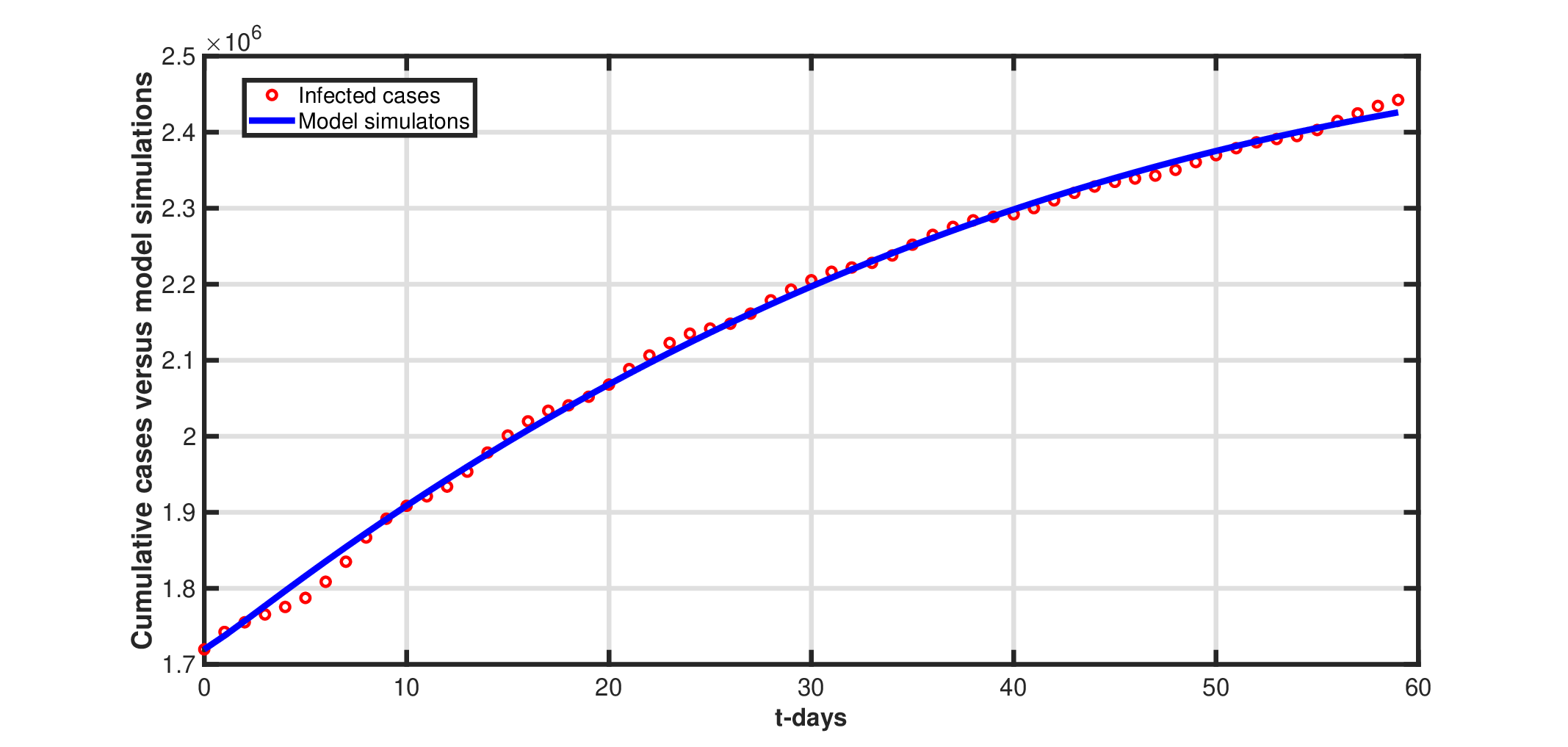}
\caption{Cumulative reported cases in Germany versus model fitting. $t=0$ stands for December 31, 2020 and $t=59$ stands for  February 28, 2021. \label{Model_fitting_germany_Vaccin}}
\end{center}
\end{figure}

The initial conditions subject to data fitting are: $S(0)=83674478$, $V(0)=49939$,  $E(0)=22924$, $A(0)=22920$,  $I(0)=32552$,  and $R(0)=97660$. The nonlinear square method is used to fit the model to the real data. Provides realistic values of model parameters, which is beneficial when we want to forecast the evolution of the disease in a given time interval. We perform experiments until the desired accurate fitting of the model is achieved. After numerically solving the optimization problem
\begin{equation}
\label{optim}
\min\limits_{\Gamma}\parallel I_{\text{predict}}-I_{\text{data}}\parallel,
\end{equation}
where $\Gamma=\left\lbrace \beta,\phi_{2},r_{1},a_{1},c_2,\eta,\theta,\gamma,\right\rbrace $, we obtain the results in Table \ref{tab_param_germany_Vacc}. The model simulation versus data fitting is shown in Figure \ref{Model_fitting_germany_Vaccin}. The value of the control reproduction number computed with the parameter values in Table \ref{tab_param_germany_Vacc} is $\mathcal{R}_c=1.127472860225384$.

\begin{table}%[t]
	%		\begin{center}
		\caption{Model parameters and their fitted values using COVID-19 data of Germany.}
		\label{tab_param_germany_Vacc}
		\begin{tabular}{@{}llllll@{}}
			\hline\noalign{\smallskip}
			Parameter & Value&Source&Parameter & Value&Source \\
			& per day&& & per day&\\
			\noalign{\smallskip}\hline\noalign{\smallskip}
			$\Lambda$ & $N(0)\times\mu$&\cite{sitePopulGerman}&$\sigma$&0.1428&\cite{Barbarossa2020}\\
			$\mu$& $\dfrac{1}{81.72\times365}$&\cite{sitePopulGerman}&$\phi_1$&0.52&\cite{siteGermanVacc}\\
		$\beta$&0.92429&Fitted&$c_1$&0.77&\cite{siteWorldVacc}\\	
			$\phi_{2}$&0.00062&Fitted&$\delta$&$0.0018$&\cite{mwalili2020seir}\\
			$r_1$&0.02534&Fitted&$r_2$&$1-r_1$&From Eq.\eqref{Covid19_model-reduce}\\
			$a_{1}$&0.34949&Fitted&$a_{2}$&$1-a_1$&From Eq.\eqref{Covid19_model-reduce}\\
			$c_2$&0.18564&Fitted&$\theta$&0.557148&Fitted\\
			$\eta$&0.35625&Fitted&$\gamma$&0.01729&Fitted\\
			$p$&0.2&Assumed&$q$&$1-p$&From Eq.\eqref{Covid19_model-reduce}\\
			\hline
		\end{tabular}
		%\end{center}
	\end{table}

\begin{figure}[ht]
\begin{center}
\includegraphics[width=1\textwidth]{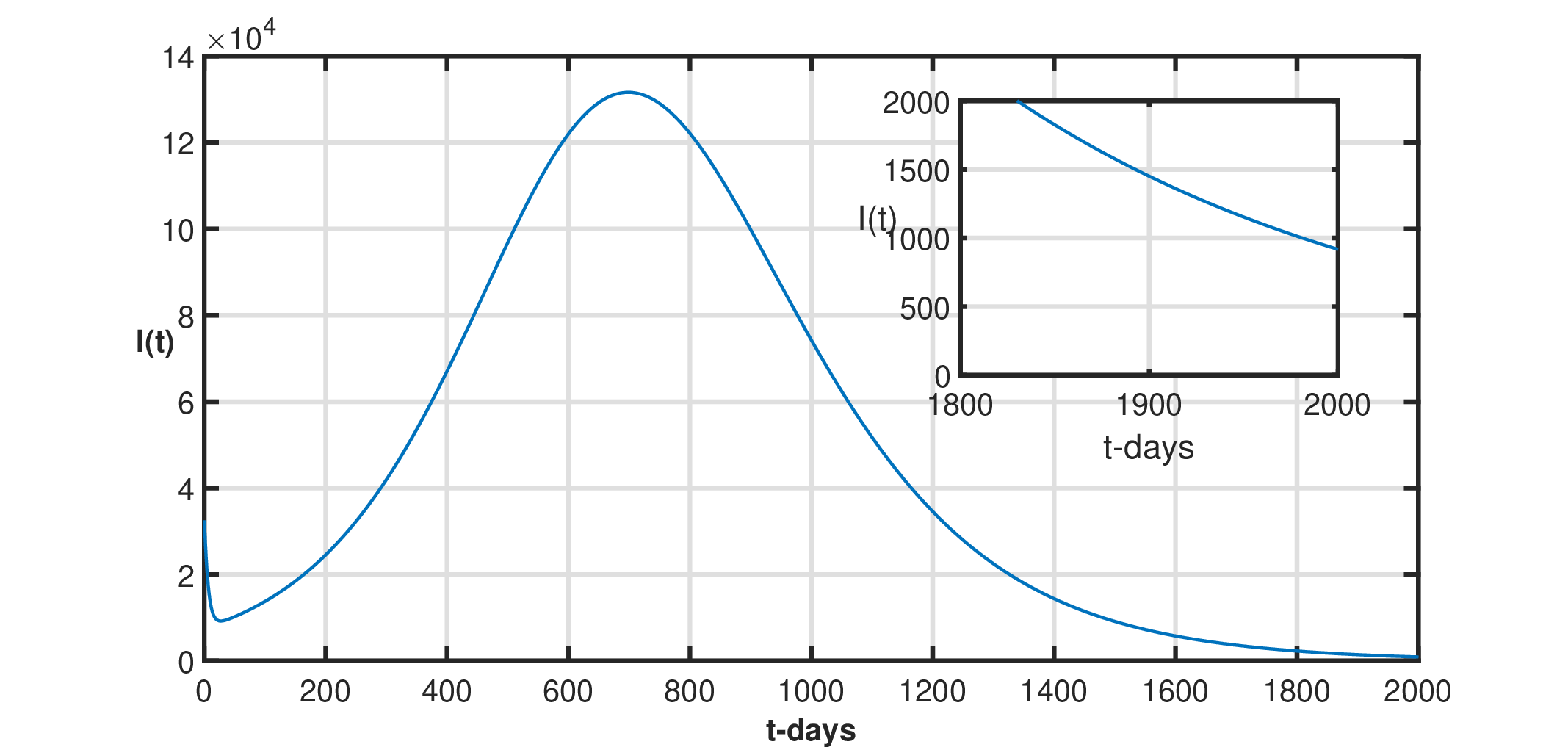}
\includegraphics[width=1\textwidth]{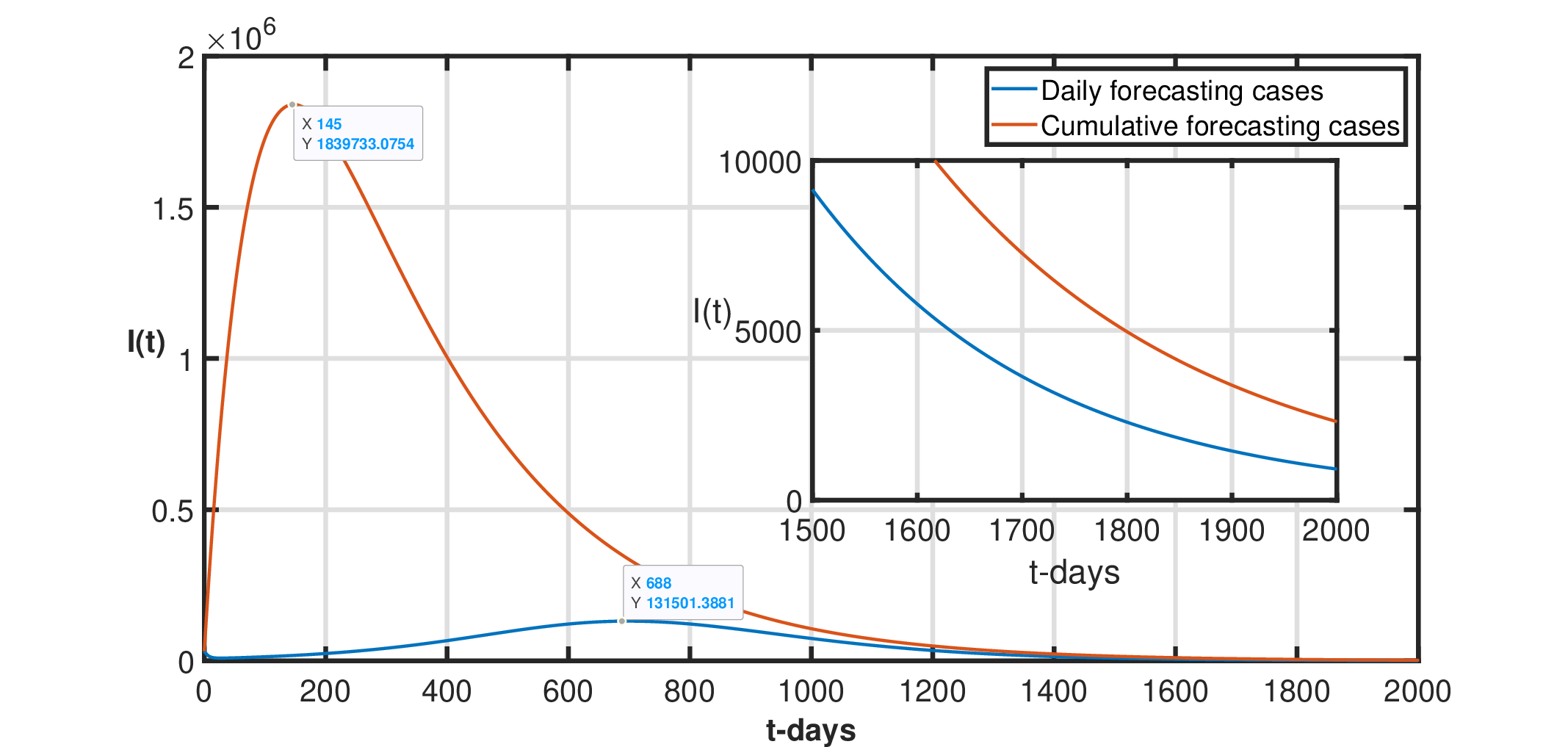}
\caption{Long-term Covid-19 prediction.\label{Prediction_germany}}
\end{center}
\end{figure}

Figure \ref{Prediction_germany} depicts the long-term prediction of COVID-19 in Germany.
\subsection{Model calibration with Cameroon data}
\begin{figure}[ht]
\begin{center}
\includegraphics[width=1\textwidth]{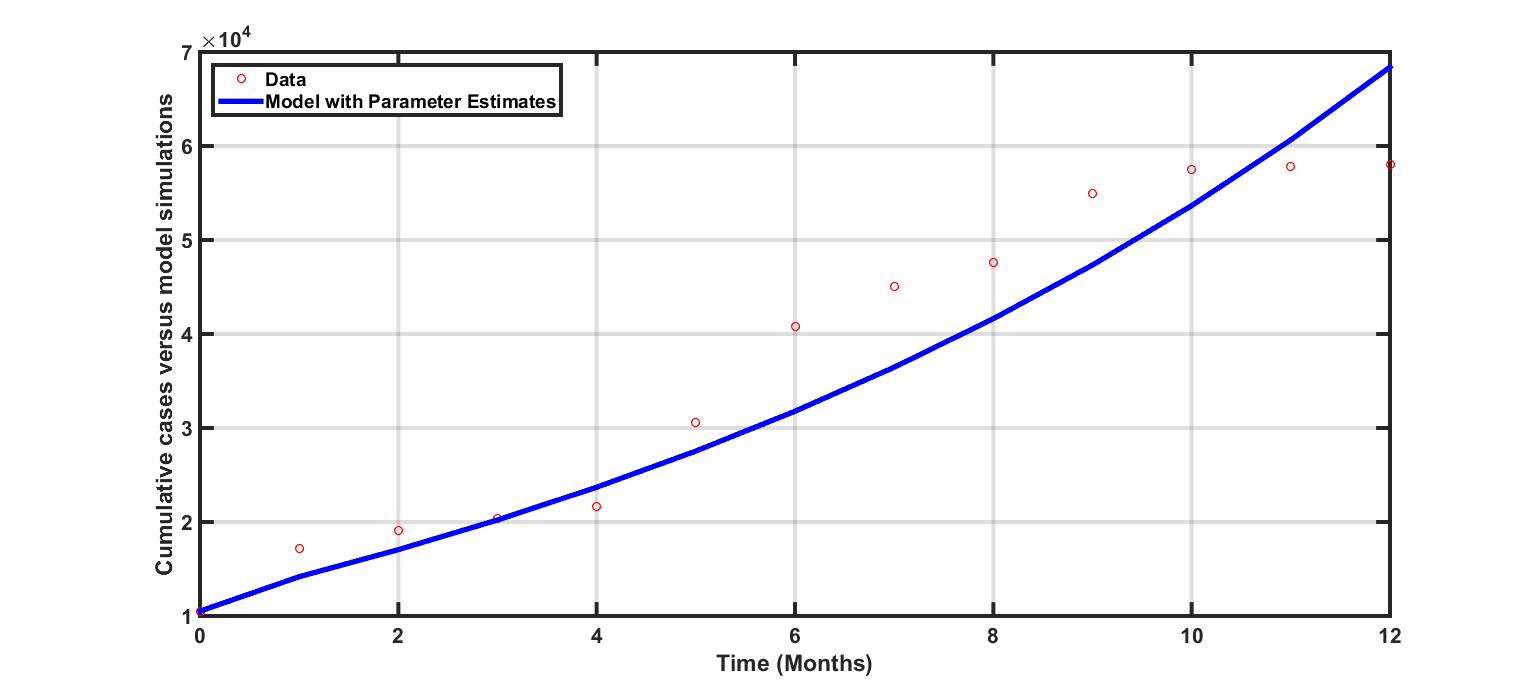}
	\caption{Cumulative reported cases in Cameroon versus model fitting. $t=0$ stands for April 2021 and $t=12$ stands for  April 2022. \label{Model_fitting_Cameroon_Vaccin}}
\end{center}
\end{figure}

In this part, we consider data from Cameroon, a Sub-Saharan country, between April 2021 and April 2022. The total number of Cameroon population is estimated to 30 millions between 2021 to 2022 \cite{Kamer_data}. The initial conditions are: $S(0)=29982802$, $V(0)=0$, $E(0)=6679$, $A(0)=0$, $I(0)=10519$, $R(0)=0$. The result is displayed in table \ref{tab_param_Cam_Vacc} and figure \ref{Model_fitting_Cameroon_Vaccin}. With these values, the control reproduction number is $\mathcal R_c \approx  1.2554$.
\begin{table}%[t]
\caption{Model parameters and their fitted values using COVID-19 data of Cameroon.}
\label{tab_param_Cam_Vacc}
		\begin{tabular}{@{}llllll@{}}
			\hline\noalign{\smallskip}
			Parameter & Value&Source&Parameter & Value&Source \\
			& per day&& & per day&\\
			\noalign{\smallskip}\hline\noalign{\smallskip}
			$\Lambda$ & $N(0)\times\mu$&\cite{Population_Cameroun}&$\sigma$&0.1428&\cite{Barbarossa2020}\\
			$\mu$& $\dfrac{1}{60\times12}$&\cite{esperance_Cameroun}&$\phi_1$&0.52&\cite{siteGermanVacc}\\
$\beta$&0.399092568990682&Fitted&$c_1$&0.77&\cite{siteWorldVacc}\\	
$\phi_{2}$&0.000583942451446   &Fitted&$\delta$&$0.0018$&\cite{mwalili2020seir}\\
$r_1$&0.061216305968569 &Fitted&$r_2$&$1-r_1$&From Eq.\eqref{Covid19_model-reduce}\\
$a_{1}$&0.000000000000092   &Fitted&$a_{2}$&$1-a_1$&From Eq.\eqref{Covid19_model-reduce}\\
$c_2$&0.000000000008747&Fitted&$\theta$&0.386526046062348&Fitted\\
$\eta$&0.004502832608954&Fitted&$\gamma$&0.869896361913556&Fitted\\
$p$&0.8262&\cite{Nkamba2021}&$q$&$1-p$&From Eq.\eqref{Covid19_model-reduce}\\
\hline
\end{tabular}
%\end{center}
	\end{table}

\section{Reaction-diffusion model and its analysis}
\label{Reaction-Diffusion}
Now, we extend the ODE model \eqref{Covid19_model-reduce} by introducing a diffusion process to obtain a reaction-diffusion model. Thus, the state variables are dependent on time and space. The Covid-19 reaction-diffusion model is then given as follows:
%\subsection{Model formulation}
{\footnotesize
	\begin{equation}
		\label{Covid19_PDE_model}
		\left\lbrace \begin{array}{l}
			\left.\begin{array}{l}
				\dfrac{\partial \mathfrak{S}(x,t)}{\partial t}-\kappa_s\Delta \mathfrak{S}(x,t)=r_1\Lambda+c_2\mathfrak{V}(x,t)-\left[k_1+\lambda(x,t)
				\right] \mathfrak{S}(x,t),\\
				\dfrac{\partial \mathfrak{V}(x,t)}{\partial t}-\kappa_v\Delta \mathfrak{V}(x,t)=r_2\Lambda+c_1\mathfrak{S}(x,t)-\left[k_2+\phi_1\lambda(x,t) \right] \mathfrak{V}(x,t),\\
				\dfrac{\partial \mathfrak{E}(x,t)}{\partial t}-\kappa_e\Delta \mathfrak{E}(x,t)=\lambda(x,t)\left(\mathfrak{S}(x,t)+\phi_1\mathfrak{V}(x,t)+\phi_2\mathfrak{R}(x,t)\right)-k_3\mathfrak{E}(x,t),\\
				\dfrac{\partial \mathfrak{A}(x,t)}{\partial t}-\kappa_a\Delta \mathfrak{A}(x,t)=p\gamma \mathfrak{E}(x,t)-k_4\mathfrak{A}(x,t),\\
				\dfrac{\partial \mathfrak{I}(x,t)}{\partial t}-\kappa_i\Delta \mathfrak{I}(x,t)=q\gamma \mathfrak{E}(x,t)+a_2\sigma \mathfrak{A}(x,t)-k_5\mathfrak{I}(x,t),\\
				\dfrac{\partial \mathfrak{R}(x,t)}{\partial t}-\kappa_r\Delta \mathfrak{R}(x,t)=a_1\sigma \mathfrak{A}(x,t)+\theta \mathfrak{I}(x,t)-\left[ \mu+\phi_2\lambda(x,t)\right]\mathfrak{R}(x,t),
			\end{array}\right\rbrace x\in\Sigma\\
			\dfrac{\partial \mathfrak{S}(x,t)}{\partial n}=\dfrac{\partial \mathfrak{V}(x,t)}{\partial n}=\dfrac{\partial \mathfrak{E}(x,t)}{\partial n}=\dfrac{\partial \mathfrak{A}(x,t)}{\partial n}=\dfrac{\partial \mathfrak{I}(x,t)}{\partial n}=\dfrac{\partial \mathfrak{R}(x,t)}{\partial n}=0,\,\, x\in\partial\Sigma,
		\end{array} \right.
	\end{equation}
}
where $\lambda(x,t)=\dfrac{\beta\left(\mathfrak{A}(x,t)+\eta \mathfrak{I}(x,t)\right)}{N(t)}$,
with the following initial conditions
\begin{equation}
	\left\lbrace \label{Init_cond_pde}
	\begin{array}{ll}
		\mathfrak{S}(x,0)=\mathfrak{S}_{0}(x)>0,\,\mathfrak{V}(x,0)=\mathfrak{V}_{0}(x)\geq 0,\, \mathfrak{E}(x,0)=\mathfrak{E}_{0}(x)\geq 0,&\\
		& x\in\Sigma.\\
		\mathfrak{A}(x,0)=\mathfrak{A}_{0}(x)\geq 0,\,\mathfrak{I}(x,0)=\mathfrak{I}_{0}(x)\geq 0,\,\mathfrak{R}(x,0)=\mathfrak{R}_{0}(x)\geq 0.&
	\end{array} \right.
\end{equation}
In system \eqref{Covid19_PDE_model}, $t$ denotes a nonnegative time;
$\Sigma$ denotes a bounded domain of $\mathbb{R}^{m}$ ($m\geq 1$) with smooth boundary $\partial\Sigma$, $n$ is the outward normal to $\partial\Sigma$; $\kappa_s>0$, $\kappa_v>0$, $\kappa_e>0$, $\kappa_a>0$, $\kappa_i>0$, $\kappa_r>0$ are the diffusion rates; $X_j(x,t)$, $j\in\left\lbrace 1,2,...,6\right\rbrace$, denote the number of population in the compartment $X_j$ in position $x$ at time $t$, with $X(x,t)=\left(\mathfrak{S}(x,t),\mathfrak{V}(x,t),\mathfrak{E}(x,t),\mathfrak{A}(x,t),\mathfrak{I}(x,t),\mathfrak{R}(x,t)\right)$.

Note that in the first open quadrant, the function $f(X)=\mathfrak{S}(\mathfrak{A}+\eta \mathfrak{I})/N$ is continuous Lipschitz for each $X_j$, $j\in\left\lbrace 1,2,...,6\right\rbrace$.%, with $X=\left(\mathfrak{S},\mathfrak{V},\mathfrak{E},\mathfrak{A},\mathfrak{I},\mathfrak{R}\right)$.
Thus, we can extend its definition to the entire first quadrant by setting $f(0,X_{I})=f(X_{S},0)=0$, where $X_{S}=\left(\mathfrak{S},\mathfrak{V},\mathfrak{R}\right)$ and $X_{I}=\left(\mathfrak{E},\mathfrak{A},\mathfrak{I} \right)$.

\subsection{Qualitative analysis of the model}
We examine the basic features of the Initial Bounded Value Problem (IVBP) \eqref{Covid19_PDE_model}-\eqref{Init_cond_pde}. Since the model deals with the population of humans, all its state variables should be positive for all $t>0$. In what follows, we thus establish the existence, uniqueness, positivity, and boundedness of the solution of the model \eqref{Covid19_PDE_model}-\eqref{Init_cond_pde}.
\begin{theorem}
	\label{exist_uniq_bound}
	For any given initial conditions that satisfy \eqref{Init_cond_pde}, the solution of model \eqref{Covid19_PDE_model}-\eqref{Init_cond_pde} is nonnegative, unique, and bounded in $[0, \infty)$.
\end{theorem}
\begin{proof}
	The IVBP model \eqref{Covid19_PDE_model}-\eqref{Init_cond_pde} can be rewritten in the Banach space $\mathcal{B}=\mathit{C}\left(\overline{\Sigma}\right)$ as follows:
	\begin{equation}
		\label{Model_compact_pde}
		\left\lbrace
		\begin{array}{l}
			\dfrac{dX(t)}{dt}=\mathbf{A}X(t)+\mathbf{f}(X(t)),\quad t>0\\
			X(0)=X_0\geq \mathbf{0}_{\mathbb{R}^{6}},
		\end{array}
		\right.
	\end{equation}
	where %$X=\left(\mathfrak{S},\mathfrak{V},\mathfrak{E},\mathfrak{A},\mathfrak{I},\mathfrak{R}\right)'$,
$X_0=\left(\mathfrak{S}_0,\mathfrak{V}_0,\mathfrak{E}_0,\mathfrak{A}_0,\mathfrak{I}_0,\mathfrak{R}_0\right)'$, $\mathbf{f}=\left(\mathbf{f}_{1},\mathbf{f}_{2},\mathbf{f}_{3},\mathbf{f}_{4},\mathbf{f}_{5},\mathbf{f}_{6}\right)'$ and formally\\
	$\mathbf{A}X=\left(\kappa_s\Delta \mathfrak{S},\kappa_v\Delta \mathfrak{V},\kappa_e\Delta \mathfrak{E},\kappa_a\Delta \mathfrak{A}, \kappa_i\Delta \mathfrak{I},\kappa_r\Delta \mathfrak{R}\right)'$, with
	\[
	\left\lbrace \begin{array}{lll}
		\mathbf{f}_{1}(X)&=&r_1\Lambda+c_2\mathfrak{V}(x,t)-\left[k_1
		+\lambda(x,t)\right] \mathfrak{S}(x,t),\\
		\mathbf{f}_{2}(X)&=&r_2\Lambda+c_1\mathfrak{S}(x,t)-\left[k_2+\phi_1\lambda(x,t) \right] \mathfrak{V}(x,t),\\
		\mathbf{f}_{3}(X)&=&\lambda(x,t)\left(\mathfrak{S}(x,t)+\phi_1\mathfrak{V}(x,t)+\phi_2\mathfrak{R}(x,t)\right) -k_3\mathfrak{E}(x,t),\\
		\mathbf{f}_{4}(X)&=&p\gamma \mathfrak{E}(x,t)-k_4\mathfrak{A}(x,t),\\
		\mathbf{f}_{5}(X)&=&q\gamma \mathfrak{E}(x,t)+a_2\sigma \mathfrak{A}(x,t)-k_5\mathfrak{I}(x,t),\\
		\mathbf{f}_{6}(X)&=&a_1\sigma \mathfrak{A}(x,t)+\theta \mathfrak{I}(x,t)-\left[ \mu+\phi_2\lambda(x,t)\right]\mathfrak{R}(x,t).
	\end{array} \right.
	\]
	$\mathbf{f}$ is locally Lipschitz in $\mathcal{B}$.
	For the existence of a local smooth solution, see \cite[Theorem B.17]{capasso1993mathematical} or \cite{mora1983semilinear,smollerbook1983}. The positivity follows from
    \cite[Theorem 14.14]{smollerbook1983}.
	
To prove the boundedness of states variables, we will follow \cite[Proposition 2.1]{han2020qualitative}. Setting $N(x,t)=\mathfrak{S}(x,t)+\mathfrak{V}(x,t)+\mathfrak{E}(x,t)+\mathfrak{A}(x,t)+\mathfrak{I}(x,t)+\mathfrak{R}(x,t)$, we obtain
	\begin{equation}
		\begin{split}
			&\dfrac{\partial N(x,t)}{\partial t}-\kappa_{s}\Delta \mathfrak{S}(x,t)-\kappa_{v}\Delta \mathfrak{V}(x,t) -\kappa_{e}\Delta \mathfrak{E}(x,t)-\kappa_{a}\Delta \mathfrak{A}(x,t)\\
			&-\kappa_{i}\Delta \mathfrak{I}(x,t)-\kappa_{r}\Delta \mathfrak{R}(x,t)\\
			&=\Lambda-\mu N(x,t)-\delta \mathfrak{I}(x,t),\quad x\in\Sigma,\,\, t>0.
		\end{split}
	\end{equation}
	Integrating the above equation over $\Sigma$, and using the Neumann boundary condition, we obtain for  $t > 0$
	\begin{equation}
		\begin{split}
			\dfrac{d}{dt}\int_{\Sigma} N(x,t)dx&=\Lambda\int_{\Sigma}dx-\mu\int_{\Sigma} N(x,t)dx-\delta \int_{\Sigma}I(x,t)dx\\
			&\leq\Lambda\left|\Sigma\right| -\mu\int_{\Sigma} N(x,t)dx
		\end{split}
	\end{equation}
	Integrating the above inequality from 0 and $t>0$ yields
	\begin{equation}
		\label{ineq_pde_func}
		\begin{split}
			\int_{\Sigma}N(x,t)dx&\leq\exp\left(-\mu t\right)\int_{\Sigma}N_{0}dx+\dfrac{\Lambda\left|\Sigma\right|}{\mu}\left[ 1-\exp(-\mu t)\right]
		\end{split}
	\end{equation}
	This implies for $j=1, 2, \ldots, 6$, that
	\begin{equation}
		\begin{split}
			\int_{\Sigma}X_{j}(x,t)dx&\leq\int_{\Sigma}N_{0}dx+\dfrac{\Lambda\left|\Sigma\right|}{\mu}.
		\end{split}
	\end{equation}
%where $X(x,t)=\left(S(x,t),V(x,t),E(x,t),A(x,t),I(x,t),R(x,t)\right)$.
Thus applying \cite[Theorem 1]{dung1997dissipativity} with $\sigma=p_{0}=1$ to the compact system \eqref{Model_compact_pde}, we conclude that there exists a positive constant $\mathcal{M}_{1}$ depending on initial data such that the solution
	$X=(S,V,E,A,I,R)$ of \eqref{Model_compact_pde} satisfies
	\begin{equation}
		\label{sup_bound_init_cond}
		\sum\limits_{j=1}^{6}\|X_{j}(.,t)\|_{L^{\infty}\left(\Sigma\right)}\leq \mathcal{M}_{1},\quad \forall t\geq 0.
	\end{equation}
	Note that \eqref{ineq_pde_func} implies
	\begin{equation}
		\begin{split}
			\lim\limits_{t\longrightarrow\infty}\int_{\Sigma}N(x,t)dx&\leq\dfrac{\Lambda\left|\Sigma\right|}{\mu}
		\end{split}
	\end{equation}
	which implies that
	\begin{equation}
		\begin{split}
			\lim\limits_{t\longrightarrow\infty}\int_{\Sigma}X_{j}(x,t)dx&\leq\dfrac{\Lambda\left|\Sigma\right|}{\mu},\quad j=1, 2, \ldots, 6.
		\end{split}
	\end{equation}
	Applying again \cite[Lemma 3.1]{peng2012reaction}, we can claim that there exists a positive constant $\mathcal{M}_{2}$ which does not depend on initial data $X^{0}_{j}$, $j=1,2,...6$ such that
	\begin{equation}
		\label{sup-bound-not-init-cond}
		\sum\limits_{j=1}^{6}\| X_{j}(.,t)\|_{L^{\infty}\left(\Sigma\right)}\leq \mathcal{M}_{2},\quad \forall t\geq T
	\end{equation}
	for some large $T>0$. This proves that each variable state of PDE model is bounded and hence, the solution exists globally, cf. \cite{capasso1993mathematical}. This ends the proof.
\end{proof}
%\begin{remark}
%Note that the case $\kappa_{s}=\kappa_{v}=...=\kappa_{r}$ of the above proof of Theorem \ref{exist_uniq_bound} can be solve as in \cite{pazzy1983semigroups} (see also Foko and Tadmon \cite{foko2022consistent}).
%\end{remark}
From Theorem \ref{exist_uniq_bound}, it follows that the following set
\[
\mathbf{W}=\left\lbrace X:=\left( \mathfrak{S},\mathfrak{V},\mathfrak{E},\mathfrak{A},\mathfrak{I},\mathfrak{R}\right)\in\mathbb{R}^{6}_{+}:0<N:=\sum\limits_{i=1}^{6}X_i\leq\frac{\Lambda\left|\Sigma\right| }{\mu}  \right\rbrace
\]
is positively invariant for system \eqref{Covid19_PDE_model}, and the IBVP \eqref{Covid19_PDE_model}-\eqref{Init_cond_pde} defines a dynamical system in it.

\subsection{Control reproduction number and stability of the DFE}
\subsubsection{The control reproduction number}
As for the ODE Covid-19 model \eqref{Covid19_model-reduce}, the PDE model \eqref{Covid19_PDE_model} always has the disease-free equilibrium \\
$\mathcal E_0=\left(\mathfrak{S}_0,\mathfrak{V}_{0},0,0,0,0\right)'$ with $\mathfrak{S}_0=\dfrac{\left(c_{2}r_{2}+r_{1}k_{2}\right)\Lambda}{\mu\left(\mu+c_{2}+c_{1}\right)}$ and $\mathfrak{V}_0=\dfrac{\left(k_{1}r_{2}+c_{1}r_{1}\right)\Lambda}{\mu\left(\mu+c_{2}+c_{1}\right)}$.
The control reproduction number is given from \eqref{Rc} by:
\begin{equation}
	\begin{split}
		\label{Rc_pde}
		\mathcal{R}_{c}&:=\dfrac{N_{1}\beta\eta\gamma\left(a_{2}p\sigma+k_{4}q\right)}{N_{0}k_{3}k_{4}k_{5}}+\dfrac{N_{1}\beta p\gamma}{N_{0}k_{3}k_{4}},
	\end{split}
\end{equation}
where $N_1=\mathfrak{S}_0+\phi_1\mathfrak{V}_0$ and $N_0=\mathfrak{S}_0+\mathfrak{V}_0=\dfrac{\Lambda}{\mu}$.
\subsubsection{Local stability of the disease-free equilibrium point}
Setting $X_1=\left(\mathfrak{S}_{1},\mathfrak{V}_{1},\mathfrak{E}_{1},\mathfrak{A}_{1},\mathfrak{I}_{1},\mathfrak{R}_{1}\right)$ where $\mathfrak{S}_{1}=\mathfrak{S}-\mathfrak{S}_{0}$, $\mathfrak{V}_{1}=\mathfrak{V}-\mathfrak{V}_{0}$, $\mathfrak{E}_{1}=\mathfrak{E}$, $\mathfrak{A}_{1}=\mathfrak{A}$, $\mathfrak{I}_{1}=\mathfrak{I}$, and $\mathfrak{R}_{1}=\mathfrak{R}$, system \eqref{Covid19_PDE_model} can be linearized as follows:
\begin{equation}
	\label{pde_linearized}
	\dfrac{\partial X_1}{\partial t}(x,t)=\mathbb{L}\left( X_{1}(x,t)\right) =\left( \varUpsilon\Delta+\mathcal{J}\left(\mathcal{E}_0\right)\right) X_{1}(x,t),
\end{equation}
where $\varUpsilon=diag\left(\kappa_{s},\kappa_{v},\kappa_{e},\kappa_{a},\kappa_{i},\kappa_{r}\right)$ and\\
$
\mathcal{J}\left(\mathcal{E}_0\right)=\left(\begin{array}{cccccc}
	-k_{1}&c_{2}&0&-\beta\dfrac{\mathfrak{S}_{0}}{N_{0}}&-\beta\eta\dfrac{\mathfrak{S}_{0}}{N_{0}}&0\\
	c_{1}&-k_{2}&0&-\beta\dfrac{\phi \mathfrak{V}_{0}}{N_{0}}&-\phi\beta\eta\dfrac{\phi \mathfrak{V}_{0}}{N_{0}}&0\\
	0&0&-k_{3}&\beta\dfrac{\mathfrak{S}_{0}+\phi \mathfrak{V}_{0}}{N_{0}}&\beta\eta\dfrac{\mathfrak{S}_{0}+\phi \mathfrak{V}_{0}}{N_{0}}&0\\
	0&0&p\gamma&-k_{4}&0&0\\
	0&0&q\gamma&a_{2}\sigma&-k_{5}&0\\
	0&0&0&a_{1}\sigma&\theta&-\mu
\end{array}\right)
$.

Let us denote by $\mathcal{I}_{6}$ the sixth-order identity matrix, $0=\psi_1\leq\psi_{2}\leq \psi_3....$, the eigenvalues of $-\Delta$ on $\overline{\mathbf{W}}$ with the zero boundary condition, and $\nu$ an eigenvalue of $\mathbb{L}$. Then, for all $i\in\mathbb{N}$, the $i$-characteristic polynomial of the operator $\mathbb{L}$ is given by
{\small
	\begin{equation}
		\label{char-Pol}
		\begin{split}
			&\forall i\in\mathbb{N},\,\,\,\\
			&\mathcal{P}_{\mathbb{L}}\left(\nu\right):=\det\left[ \mathcal{J}\left(\mathcal{E}_0\right)-\psi_{i}\varUpsilon-\nu\mathcal{I}_{6}\right]\\
&=\left(\nu+\psi{i}\kappa_{r}+\mu\right)\\
&\times\overbrace{\left[\nu^2+\left(\psi_{i}(k_{v}+k_{s})+k_{2}+k_{1}\right)\nu+
	(\left(\psi_{i}k_{s}+k_{1}\right)k_{v}+k_{2}k_{s})\psi_{i}+\mu\left(\mu+c_{1}+c_{2}\right) \right]}\limits^{g^{(1)}_{i}(\nu)}\\
&\quad\times
\overbrace{\left[
	\nu^{3}+\varpi_{2}x^{2}+\varpi_{1}\nu+\varpi_{0}
	\right]}\limits^{g^{(2)}_{i}(\nu)}
\end{split}
\end{equation}
}
where
$\varpi_{2}=\left(\kappa_{i}+\kappa_{e}+\kappa_{a}\right)\psi_{i}+k_{5}+k_{4}+k_{3}$,
%\footnotesize
\begin{equation*}
\begin{split}
\varpi_{1}&=\frac{1}{a_{2}\eta p\sigma+k_{4}\eta q+k_{5}p}\times\\
&\left\lbrace\left[ \left(\left(\kappa_{e}+\kappa_{a}\right) \kappa_{i}+\kappa_{a}\kappa_{e}\right)\psi_{i}
+\left(k_{4}+k_{3}\right)\kappa_{i
}+\left(k_{5}+k_{4}\right)\kappa_{e}
+\left(k_{5}+k_{3}\right)\kappa_{a}
\right]\right.\\
&\times \left.\left( a_{2}\eta\psi_{i} p\sigma+k_{4}\eta\psi_{i} q+k_{5}\psi_{i}p\right) \right.\\
&\left.+\left( 1-\mathcal R_{c}\right)\left(p+\eta q\right) k_{3}k_{4}k_{5}
+\left(k_{4}+k_{3}\right)\left( k_{5}+a_{2}\eta \sigma\right)k_{5} p\right.\\
&\left.
+\left(\left(k_{5}+k_{3}\right) k_{4}q+k_{3}a_{2}p\sigma\right)k_{4}\eta\right\rbrace,\\
\varpi_{0}&=\frac{1}{a_{2}\eta p
\sigma+k_{4}\eta q+k_{5}p}\times\\
&\left\lbrace \left[ \kappa_{a}\kappa_{e}\kappa_{i}\psi_{i}^2+\left(\left(k_{4}\kappa_{e}+k_{3}\kappa_{a}
\right)\kappa_{i}+k_{5}\kappa_{a}\kappa_{e}
\right)\psi_{i}+\left( k_{3}\kappa_{i}+k_{5}\kappa_{e}\right)k_{4}\right]\right.\\
&\times\left.\left( a_{2}\eta\psi_{i}p\sigma+k_{4}\eta\psi_{i}q\right) \right.\\
&\left.+\left[\kappa_{a}\kappa_{e}\kappa_{i}\psi_{i}^2
+\left(\left(k_{4}\kappa_{e}+k_{3}\kappa_{a}\right)\kappa_{i}+k_{5}\kappa_{a}\kappa_{e}\right)\psi_{i}+\left( k_{4}\kappa_{e}+k_{3}\kappa_{a}\right)k_{5}\right]\psi_{i}k_{5} p\right.\\
&\left.+k_{3}k_{5}\kappa_{a}a_{2}\eta\psi_{i}p\sigma
+\left(1-\mathcal R_{c}\right)\left[ a_{2}\eta p\sigma +\left(\kappa_{i}\psi_{i}+k_{5}\right)p+\left(\kappa_{a}\psi_{i} +k_{4}\right) \eta q\right]k_{3}k_{4}k_{5}\right\rbrace
\end{split}
\end{equation*}
It is clear that for any $i\geq 1$, the roots of $\mathcal{P}_{\mathbb{L}}$ are $\nu=-\psi{i}\kappa_{r}-\mu$, and those of $g^{(1)}_{i}$ and $g^{(2)}_{i}$. Note that $\mathcal R_c<1$ ensures that $\varpi_1$ and $\varpi_0$ are positive. Since for $i\geq 1$, $\psi_i\geq 0$, it follows that the roots of $g^{(1)}_{i}$ have negative real parts.  A rigorous algebraic computation permits to prove that $\varpi_{2}\varpi_{1}-\varpi_{0}>0$. Thus, $g_{i}^{(2)}$ meets the Routh-Hurwitz criteria. It then follows that all roots of the characteristic polynomial \eqref{char-Pol} have negative real parts whenever $\mathcal R_c<1$, which implies that the disease-free equilibrium $\mathbb E_0$ is locally stable (see \cite{tadmon2022transmission}).

Since $\psi_1=0$, it follows that, if $\mathcal{R}_c>1$,\\
$
g_{1}^{(2)}(0)=\left(1-\mathcal{R}_{c}\right)\left[ a_{2}\eta p\sigma +\left(\kappa_{i}\psi_{i}+k_{5}\right)p+\left(\kappa_{a}\psi_{i} +k_{4}\right) \eta q\right]k_{3}k_{4}k_{5}<0$\\ and $\lim\limits_{\nu\rightarrow+\infty}g_{1}^{(2)}(\nu)=+\infty$,
which implies that the characteristic polynomial $\mathcal{P}_{\mathbb{L}}$ admits at least one root with positive real part. Thus, the disease-free equilibrium $\mathbb E_0$ is not stable.
The above analysis can be summarized as follows:
\begin{theorem}
\label{las_dfe_pde_covid}
The disease-free equilibrium $\mathbb E_0$ of the PDE model \eqref{Covid19_PDE_model} is locally asymptotically stable in
$\mathbf{W}$ if $\mathcal R_c<1$, and unstable otherwise.
\end{theorem}

\subsubsection{Global stability of the disease-free equilibrium}
\begin{theorem}
If $\mathcal{R}_{c}<1$, then the disease free equilibrium  $\mathbb{E}_0$ of the COVID-19 PDE model \eqref{Covid19_PDE_model} is globally asymptotically stable provided that
\begin{equation}
\label{Cond_gas_dfe}
\dfrac{N_{1}}{N_{0}}-\dfrac{\left(\mathfrak{S}(x,t)+\phi_1\mathfrak{V}(x,t)+\phi_2\mathfrak{R}(x,t)\right)}{N(t)}\geq 0.
\end{equation}
\end{theorem}
\begin{proof}
Let us denote by $\left(\mathfrak{S}(x,t),\mathfrak{V}(x,t),\mathfrak{E}(x,t),\mathfrak{A}(x,t),\mathfrak{I}(x,t),\mathfrak{R}(x,t)\right)'$ any arbitrary nonnegative solution of COVID-19 IVBP \eqref{Covid19_PDE_model}-\eqref{Init_cond_pde}, $N_0=\mathfrak{S}_0+\mathfrak{V}_0$, $N_1=\mathfrak{S}_0+\phi_1\mathfrak{V}_0$, with $\mathfrak{S}_0=\dfrac{\left(c_{2}r_{2}+r_{1}k_{2}\right)\Lambda}{\mu\left(\mu+c_{2}+c_{1}\right)}$ and $\mathfrak{V}_0=\dfrac{\left(k_{1}r_{2}+c_{1}r_{1}\right)\Lambda}{\mu\left(\mu+c_{2}+c_{1}\right)}$.
To prove the global asymptotic stability of the disease-free equilibrium $\mathcal{E}_0=\left(\mathfrak{S}_{0},\mathfrak{V}_{0},0,0,0,0\right)'$, we consider %as for the case of the Covid-19 ode model \eqref{Covid19_model-reduce},
the following Lyapunov functional $\mathcal{K}(t)$ defined as
\begin{equation*}
\begin{split}
\mathcal{K}(\mathfrak{S},\mathfrak{V},\mathfrak{E},\mathfrak{A},\mathfrak{I},\mathfrak{R})&=\int_{\Sigma}\mathcal{Q}\left(\mathfrak{S}(x,t),\mathfrak{V}(x,t),\mathfrak{E}(x,t),\mathfrak{A}(x,t),\mathfrak{I}(x,t),\mathfrak{R}(x,t)\right) dx\\
&=\int_{\Sigma}\left(w'W^{-1}\left(\mathfrak{E},\mathfrak{A},\mathfrak{I}\right)\right) dx\\
\end{split}
\end{equation*}
where $w'$ is the left eigenvector of the nonnegative matrix $W^{-1}Z$ corresponding to the eigenvalue $\mathcal{R}_c$. From \cite[Proposition 2.1]{hattaf2013global}, the time derivative of $\mathcal{K}$ along the nonnegative solution of model \eqref{Covid19_PDE_model} is given by
\begin{equation*}
\begin{split}
\dfrac{d\mathcal{K}}{dt}(\mathfrak{S},\mathfrak{V},\mathfrak{E},\mathfrak{A},\mathfrak{I},\mathfrak{R})	
&=\int_{\Sigma}\left[\left(\mathcal{R}_c-1\right)w'\left(\mathfrak{E},\mathfrak{A},\mathfrak{I}\right)-w'W^{-1}\mathcal{M}\left(\mathfrak{S},\mathfrak{V},\mathfrak{E},\mathfrak{A},\mathfrak{I},\mathfrak{R}\right)\right] dx\\
\end{split}
\end{equation*}
This implies that if $\mathcal{R}_c<1$, $\dfrac{d\mathcal{K}}{dt}\leq 0$ whenever $w'W^{-1}\mathcal{M}\left(\mathfrak{S},\mathfrak{V},\mathfrak{E},\mathfrak{A},\mathfrak{I},\mathfrak{R}\right)\geq\mathbf{0}_{\mathbb{R}^{3}}$, which is equivalent to the condition  $\dfrac{N_{1}}{N_{0}}-\dfrac{\left(\mathfrak{S}(x,t)+\phi_1\mathfrak{V}(x,t)+\phi_2\mathfrak{R}(x,t)\right)}{N(t)}\geq 0$. Furthermore, $\dfrac{d\mathcal{K}}{dt}=0$ holds if and only if $\mathfrak{E}=\mathfrak{A}=\mathfrak{I}=0$. Then, $\left\lbrace\mathcal{E}_0\right\rbrace $ is the only compact invariant subset of $\left\lbrace\left(\mathfrak{S},\mathfrak{V},\mathfrak{E},\mathfrak{A},\mathfrak{I},\mathfrak{R}\right)'\in\mathbb{R}^{6}:\dfrac{d\mathcal{K}}{dt}=0\right\rbrace $. By LaSalle's invariance principle \cite{henry2006geometric,la1976stability}, if $\mathcal{R}_{c}<1$, then the disease-free equilibrium $\mathcal{E}_{0}$ of the IBVP \eqref{Covid19_PDE_model}-\eqref{Init_cond_pde} is globally asymptotically stable in $\Sigma$ provided that condition \eqref{Cond_gas_dfe} holds.
\end{proof}
\begin{remark}
Note that inequality \eqref{Cond_gas_dfe} is just to ensure the global asymptotic stability of the disease-free equilibrium, of both models (Eq. \eqref{Covid19_model-reduce} and Eq.\eqref{Covid19_PDE_model}). This condition does not have a biological interpretation. These kind of conditions are used in the literature (see, for example, \cite{li2016class,bowongmathematical} and the references therein).
\end{remark}
\subsection{Existence of endemic equilibrium points}
Let \\
$\mathcal{E}=\mathcal{E}(x)=\left(
\mathfrak{S}^{\star}, \mathfrak{V}^{\star},\mathfrak{E}^{\star},\mathfrak{A}^{\star},\mathfrak{I}^{\star},\mathfrak{R}^{\star}\right)'$ any spatially equilibrium of system \eqref{Covid19_PDE_model}. Thus $\mathcal{E}$ must solve the following system
\begin{equation}
\label{ee_pde}
\left\lbrace \begin{array}{ll}
r_1\Lambda+c_2\mathfrak{V}^{\star}-k_1\mathfrak{S}^{\star}-\lambda^{\star}\mathfrak{S}^{\star}&=0,\\
r_2\Lambda+c_1\mathfrak{S}^{\star}-\left[k_2+\phi_1\lambda^{\star}\right] \mathfrak{V}^{\star}&=0,\\
\lambda^{\star}\left(\mathfrak{S}^{\star}+\phi_1\mathfrak{V}^{\star}+\phi_2\mathfrak{R}^{\star}\right)-k_3\mathfrak{E}^{\star}&=0 ,\\
p\gamma \mathfrak{E}^{\star}-k_4\mathfrak{A}^{\star}&=0,\\
q\gamma \mathfrak{E}^{\star}+a_2\sigma \mathfrak{A}^{\star}-k_5\mathfrak{I}^{\star}&=0,\\
a_1\sigma \mathfrak{A}^{\star}+\theta \mathfrak{I}^{\star}-\left[ \mu+\phi_2\lambda^{\star}\right]\mathfrak{R}^{\star}&=0\\
\end{array} \right.
\end{equation}
Note that algebraic system \eqref{ee_pde} is equivalent to \eqref{ee1}. Thus from \cite{hattaf2013global}, we conclude that, as for the ODE Covid-19 model \eqref{Covid19_model-reduce}: %we have the following result:
\begin{proposition}$\;$
\begin{enumerate}
\item If $\mathcal R_c<1$, then the Covid-19 PDE model \eqref{Covid19_model-reduce} admits a disease-free equilibrium which can co-exist with two endemic equilibrium points, depending of the sign of coefficients $\mathcal A_{2}$ and $\mathcal A_{1}$ of the polynomial \eqref{Pol_ee}.
\item If $\mathcal R_c>1$, then the disease-free equilibrium point of the Covid-19 PDE model \eqref{Covid19_model-reduce} can coexist with one or three endemic equilibrium points, depending of the sign of coefficients $\mathcal A_{2}$ and $\mathcal A_{1}$ of the polynomial \eqref{Pol_ee}.
\end{enumerate}
\end{proposition}

\section{Numerical scheme and simulations}
\label{num_analysis}
In this part of the work, we will construct a numerical scheme to simulate the Coronavirus (COVID-19) spatiotemporal model \eqref{Covid19_PDE_model}.
\subsection{Description of the numerical method}
The system of partial differential equations \eqref{Covid19_PDE_model} can be numerically solved using the Partial Differential Equation Toolbox in \textsc{Matlab}${}^\copyright$ (see \cite{matlabDoc}). Among other things, it can solve systems of equations of the form $$m u_{tt} + d u_t - \nabla \cdot (c \otimes \nabla u) + au = f,$$ using a finite elements approximation. The coefficients $m$, $d$, $c$, $a$ and $f$ are allowed to depend on time, space, and $u$ itself. In our case where $u = (S,V,E,A,I,R)$, we have $m = 0$, $d = 1$, $c = \operatorname{diag}(\kappa_S, \kappa_V, \kappa_E, \kappa_A, \kappa_I, \kappa_R)$, $f = (r_1 \Lambda, r_2 \Lambda, 0,0,0,0)$ and
\[
a(u) = -\begin{pmatrix}
-(k_1 + \lambda(u)) & c_2 & 0 & 0 & 0 & 0\\
c_1 & -(k_2 + \phi_1 \lambda(u)) & 0 & 0 & 0 & 0\\
\lambda(u) & \phi_1 \lambda(u) & -k_3 & 0 & 0 & \phi_2 \lambda(u)\\
0 & 0 & p\gamma & -k_4 & 0 & 0\\
0 & 0 & q\gamma & a_2 \sigma & -k_5 & 0\\
0 & 0 & 0 & a_1 \sigma & \theta & -(\mu + \phi_2 \lambda(u))
\end{pmatrix},
\]
where $\lambda(u) = \beta(u_4 + \eta u_5)/(u_1 + u_2 + u_3 + u_4 + u_5 + u_6)$. Together with these coefficients, we have to specify boundary conditions (Neumann is the default) and initial values. For the geometry on which to solve the equations, we use map data from GADM \cite{GADM}, reduced in complexity with the help of Mapshaper \cite{mapshaper}.

\subsection{Numerical simulations}
In this section we provide numerical results of the reaction-diffusion model \eqref{Covid19_PDE_model}. We consider two approaches: \textit{(1)} all model parameters are constant and \textit{(2)} the transmission rate coefficient $\beta$ is time-space dependent. For each case, we present the initial distribution of each variable state at $t=0$, and the final distribution of each variable state of the model at $t=t_{final}$.
\subsubsection{Germany}
Here we assume that $\kappa_s=\kappa_v=\kappa_e=\kappa_a=\kappa_r=1$ and $\kappa_i=0.1$, which means that only individuals who have been tested positive are confined (or their movements are restricted). Note that we assume that the confinement or the restriction is not perfect, hence $\kappa_i\ne 0$.

When given spatially constant initial values, the PDE algorithm yields the same results as the ODE solver, which is also to be expected from theory. In the following we therefore present different initial values in order to show the differences between the two approaches. We show the distribution of the susceptible, vaccinated, and infected individuals only.

\begin{figure}[ht]
\includegraphics[width = \textwidth]{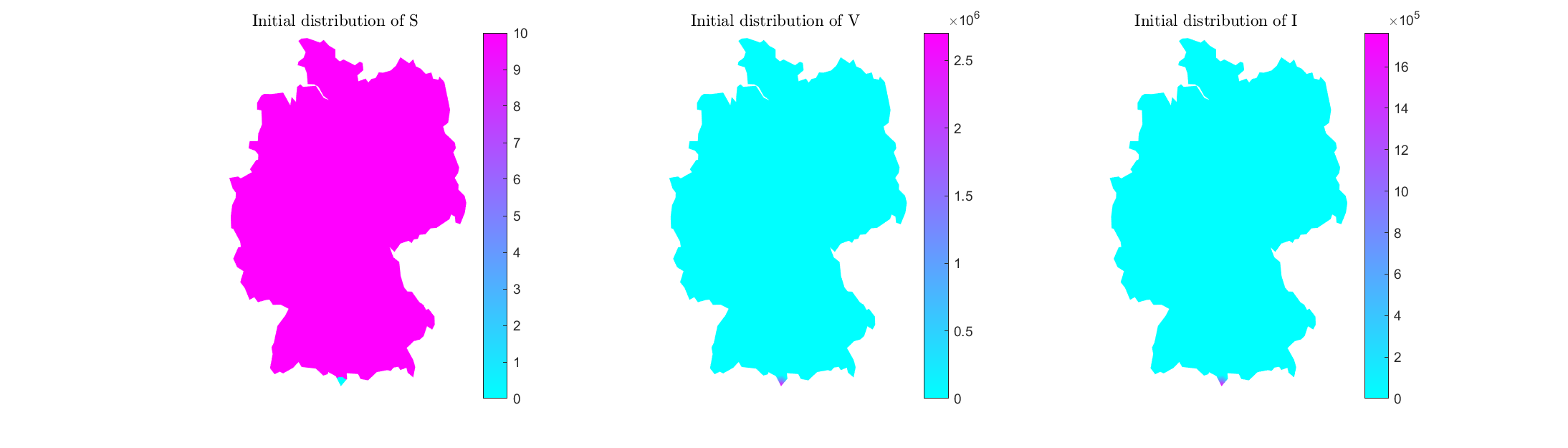}
\caption{Initial distribution of model state variables where there is one southern peak. \label{1PeaksInit}}
\end{figure}

\begin{figure}[ht]
%\begin{center}
\includegraphics[width = \textwidth]{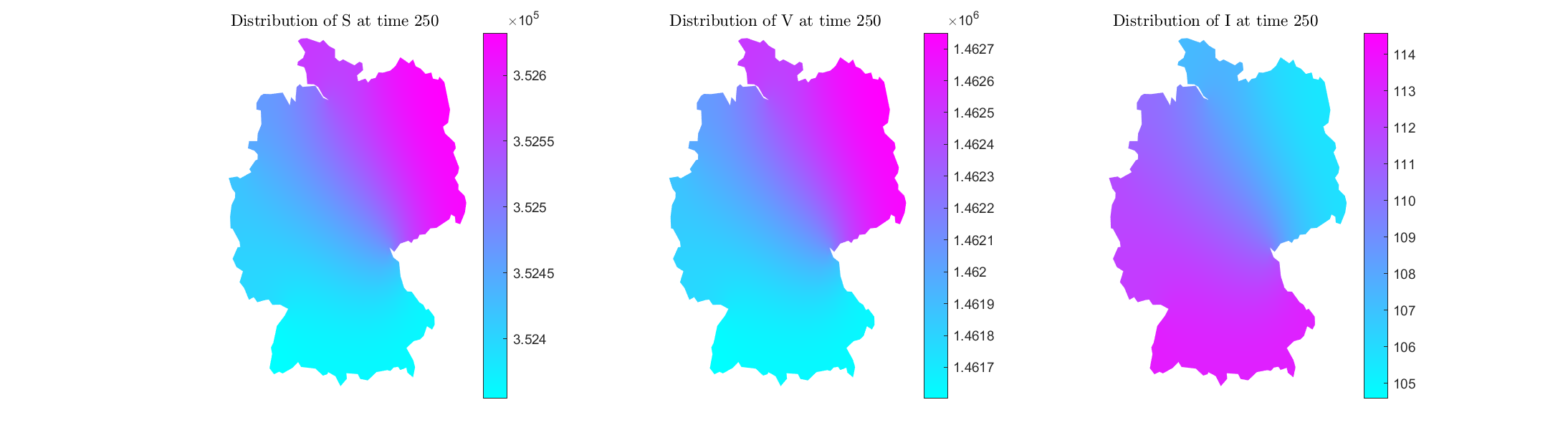}
\caption{Solution after $t_f=250$ days. \label{1peak250}}
%\end{center}
\end{figure}

\begin{figure}[ht]
%\begin{center}
\includegraphics[width = \textwidth]{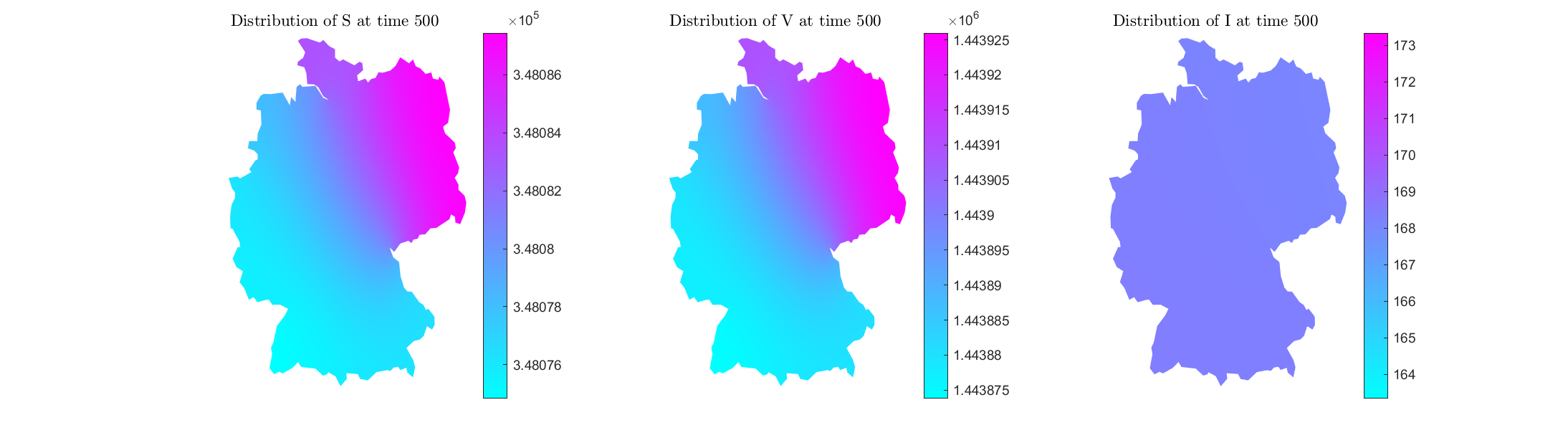}
\caption{Solution after $t_f=500$ days. \label{1peak500}}
%\end{center}
\end{figure}

\begin{figure}[ht]
\begin{center}
\includegraphics[width = \textwidth]{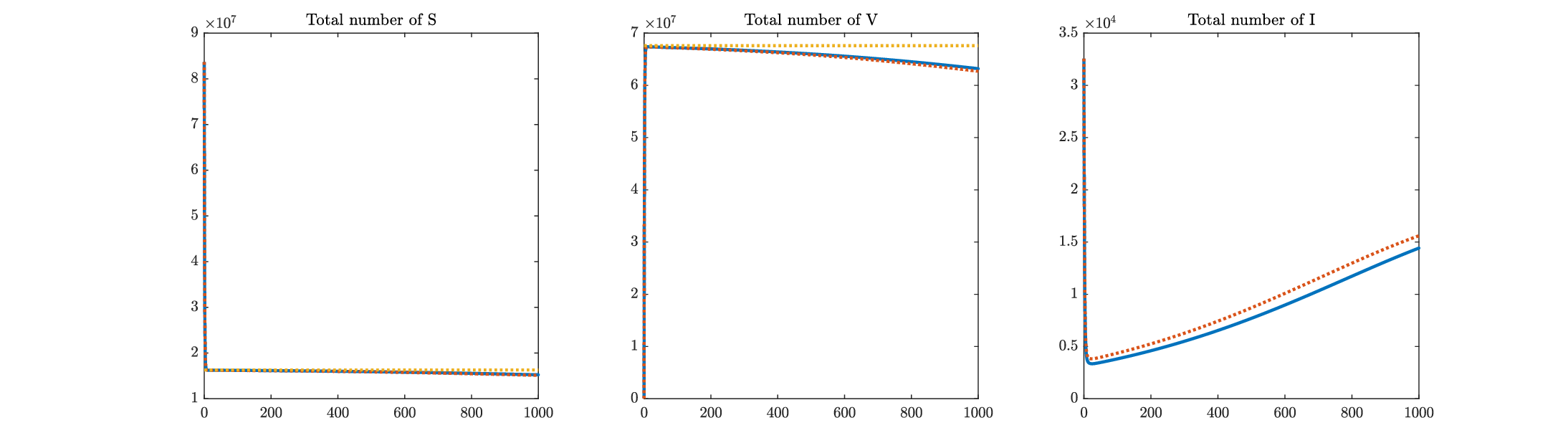}
\caption{Total number of members in each compartment over time (blue), compared to the ODE model (red) and the disease-free equilibrium (yellow). \label{1PeaksTotal}}
\end{center}
\end{figure}

\paragraph{One peak}
For the first test we chose the initial population to be all susceptible to infection, except for one small region in the very south of Germany, where there are also infected, exposed,~\ldots persons. The values are adjusted in such a way that the total number of members in every compartment is the same as before, i.\,e. $S(0)=83674478$, $V(0)=49939$,  $E(0)=22924$, $A(0)=22920$,  $I(0)=32552$,  and $R(0)=97660$, as in the beginning of Section \ref{calibration}. We use the coefficients from Table \ref{tab_param_germany_Vacc}. The results are displayed in figures \ref{1PeaksInit}-\ref{1PeaksTotal}. Figure \ref{1PeaksInit} presents the initial distribution of the model state variables while figures \ref{1peak250}-\ref{1peak500} display solution of the PDE model after $t_f=250$ days and $t_f=500$ days, respectively. The total number of members in each compartment over time (blue), compared to the ODE model (red) and the disease-free state from Equation \ref{DFE} (yellow), is displayed in figure \ref{1PeaksTotal}.

\begin{figure}[ht]
\includegraphics[width = \textwidth]{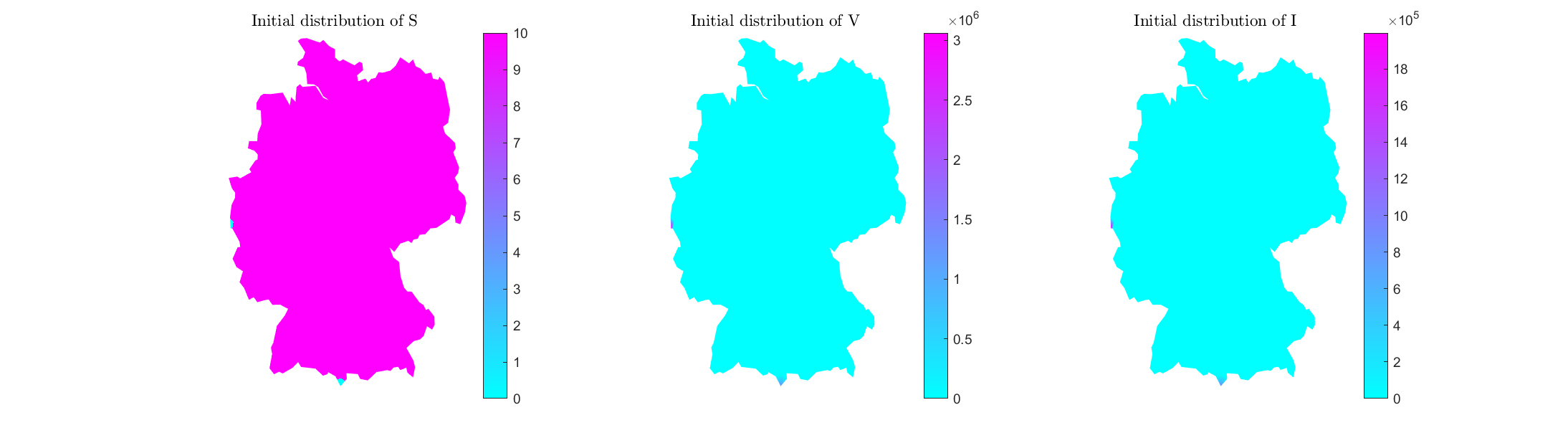}
\caption{Initial distribution of model state variables where there are two peaks (south and west). \label{2peakInit}}
\end{figure}

\begin{figure}[ht]
\includegraphics[width = \textwidth]{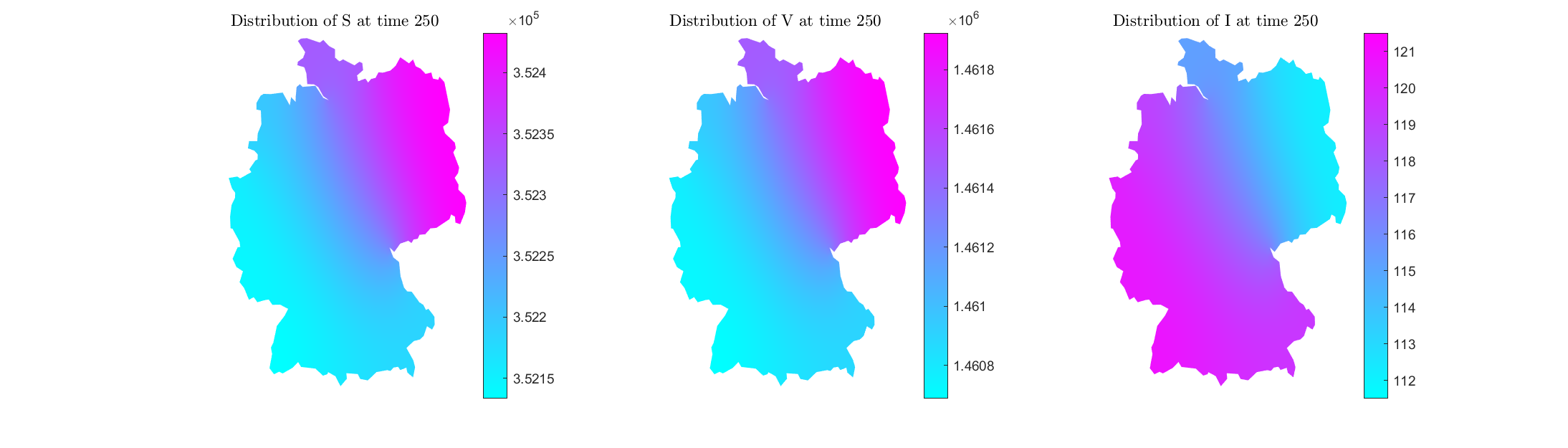}
\caption{Solution after $t_f=250$ days. \label{2peak250}}
\end{figure}

\begin{figure}[ht]
\includegraphics[width = \textwidth]{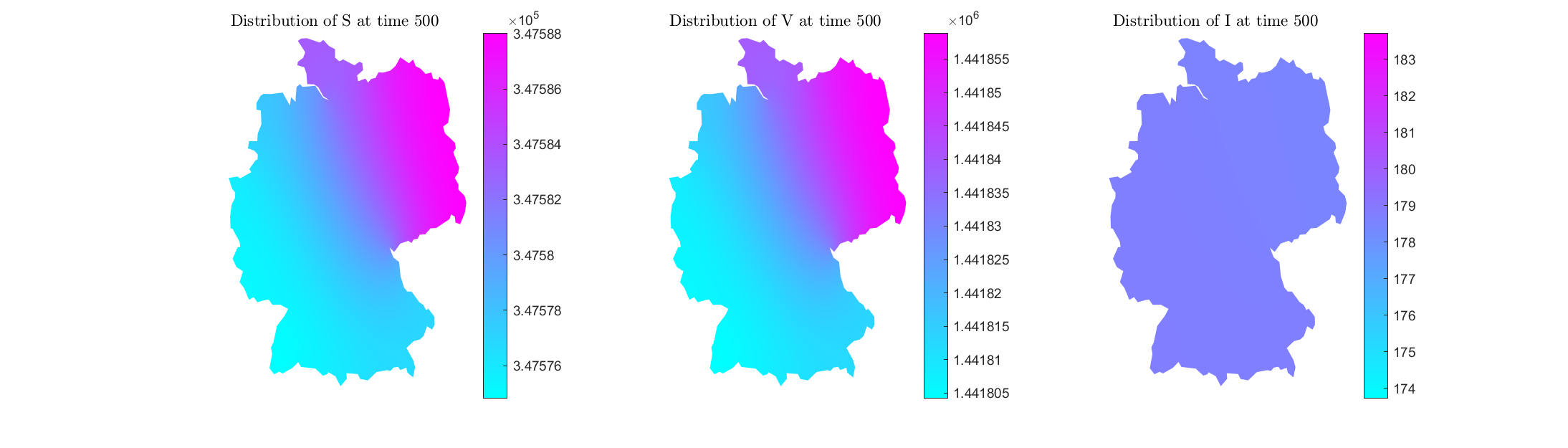}
\caption{Solution after $t_f=500$ days. \label{2peak500}}
\end{figure}

\begin{figure}[ht]
\begin{center}
\includegraphics[width = \textwidth]{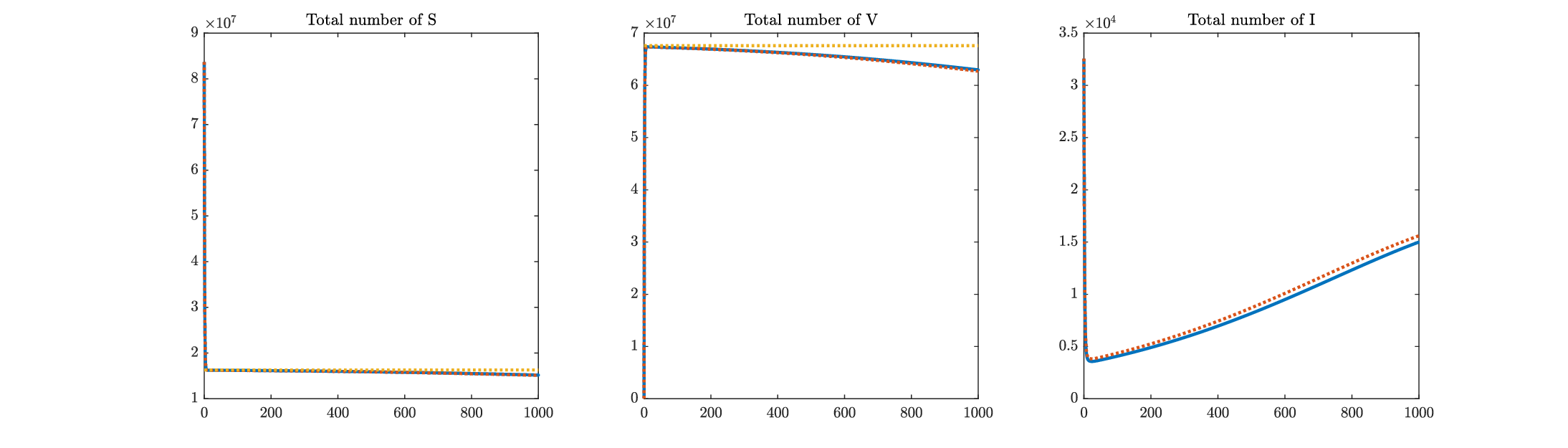}
\caption{Total number of members in each compartment over time (blue), compared to the ODE model (red) and the disease-free state (yellow). \label{2PeakTotal}}
\end{center}
\end{figure}

\begin{figure}[ht]
\includegraphics[width = \textwidth]{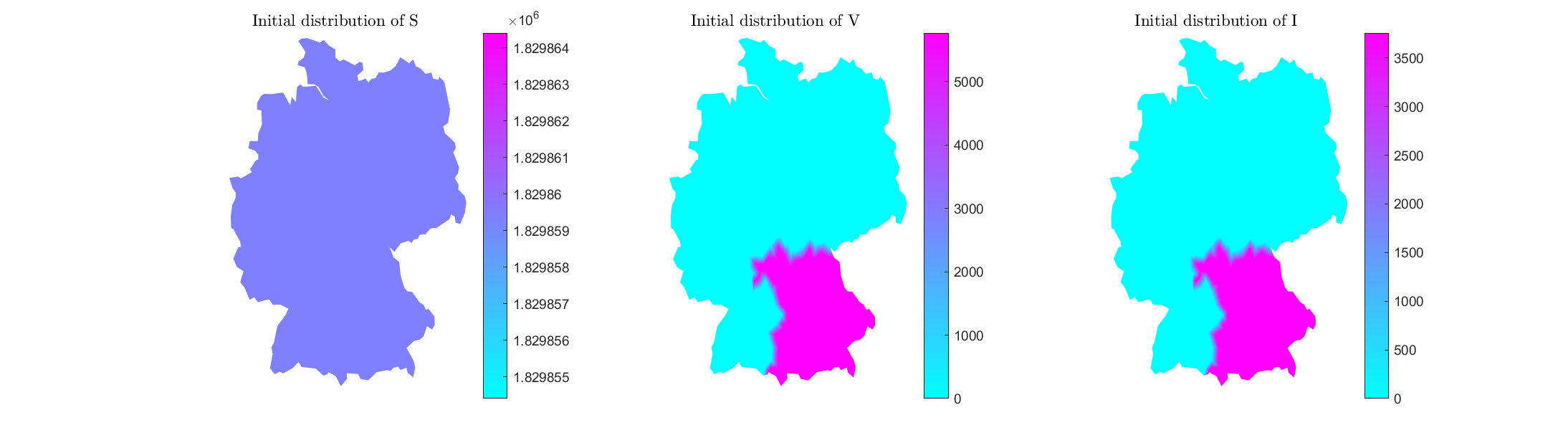}
\caption{Initial distribution of model state variables where only the state of Bavaria holds exposed, asymptomatic infected, and symptomatic infected individuals. \label{bavariaInit}}
\end{figure}

\begin{figure}[ht]
\includegraphics[width = \textwidth]{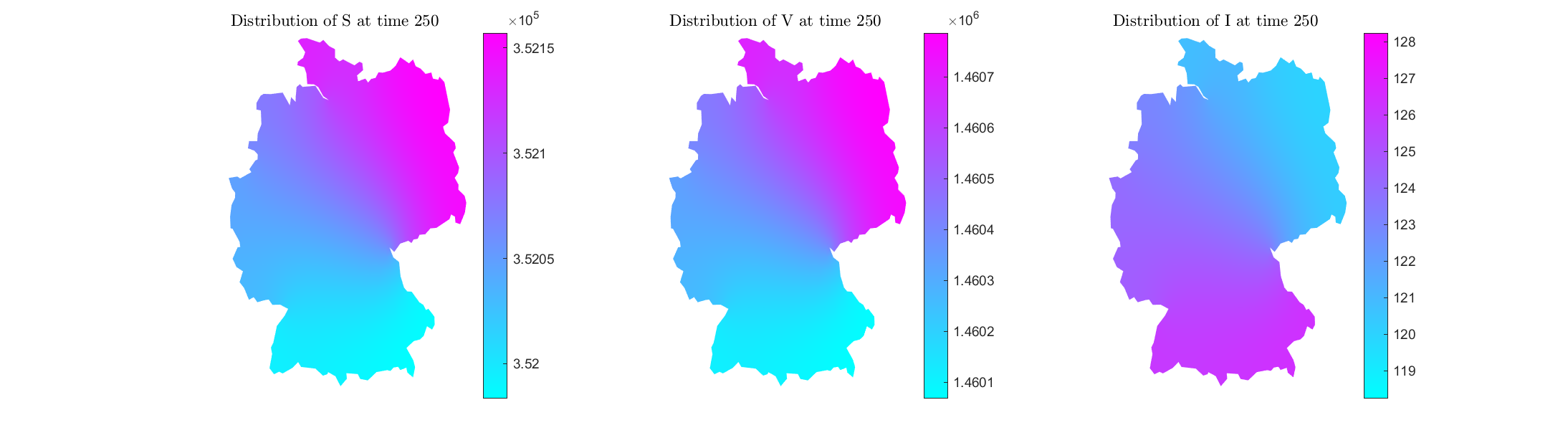}
\caption{Solution after $t_f=250$ days. \label{bavaria250}}
\end{figure}

\begin{figure}[ht]
\includegraphics[width = \textwidth]{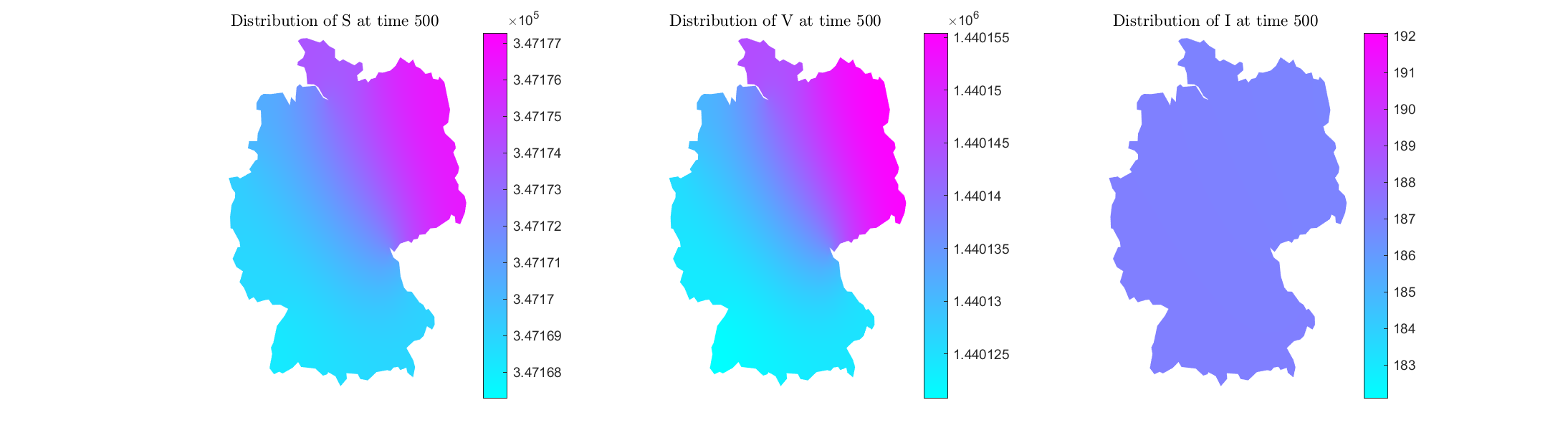}
\caption{Solution after $t_f=500$ days. \label{bavaria500}}
\end{figure}

\begin{figure}[ht]
%\begin{center}
\includegraphics[width = \textwidth]{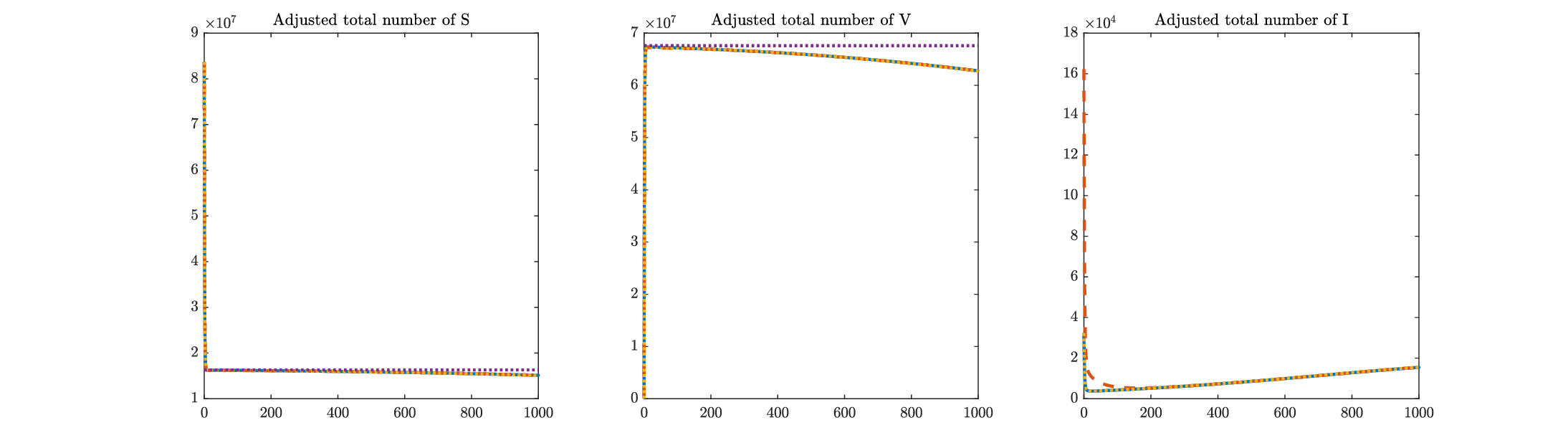}
\caption{Total number of members in each compartment over time (blue), compared to the numbers only in Bavaria adjusted for size (red), the ODE model (yellow) and the disease-free equilibrium (purple). \label{bavariaTotal}}
%\end{center}
\end{figure}

%\FloatBarrier
\paragraph{Two peaks}
For the next test we added a second peak in western Germany, close to Aachen/Heinsberg, where a major outbreak of Covid occurred in early 2020. Again, the total numbers of members in each compartment did not change. The results are display in Figures \ref{2peakInit}-\ref{2PeakTotal}. Figure \ref{2peakInit} presents the initial distribution of model state variables while Figures \ref{2peak250}-\ref{2peak500} display the solution of the PDE model after $t_f=250$ days and $t_f=500$ days, respectively. The total number of members in each compartment over time (blue), compared to the ODE model (red) and the disease-free state (yellow) is displayed in figure \ref{2PeakTotal}.
%\FloatBarrier
\paragraph{Bavaria}
Figures \ref{bavariaInit}-\ref{bavariaTotal} display a series of images where only the state of Bavaria holds exposed, asymptomatic infected, and symptomatic infected individuals. Initial distribution of model state variables are displayed on figures \ref{bavariaInit}, while  solution after $t_f=250$ days and $t_f=500$ days, respectively, are display in figure \ref{bavaria250}-\ref{bavaria500}. We see that more time passes, the more the infection spreads to the rest of the country. Note that the total numbers in each compartment are slightly different to the ones before, this is because \textsc{Matlab} smooths out the initial value along the Bavarian border. The total number of members in each compartment over time (blue), compared to the ODE model (red) and the disease-free state (yellow), is displayed in Figure \ref{bavariaTotal}.


\begin{figure}[ht]
\includegraphics[width = \textwidth]{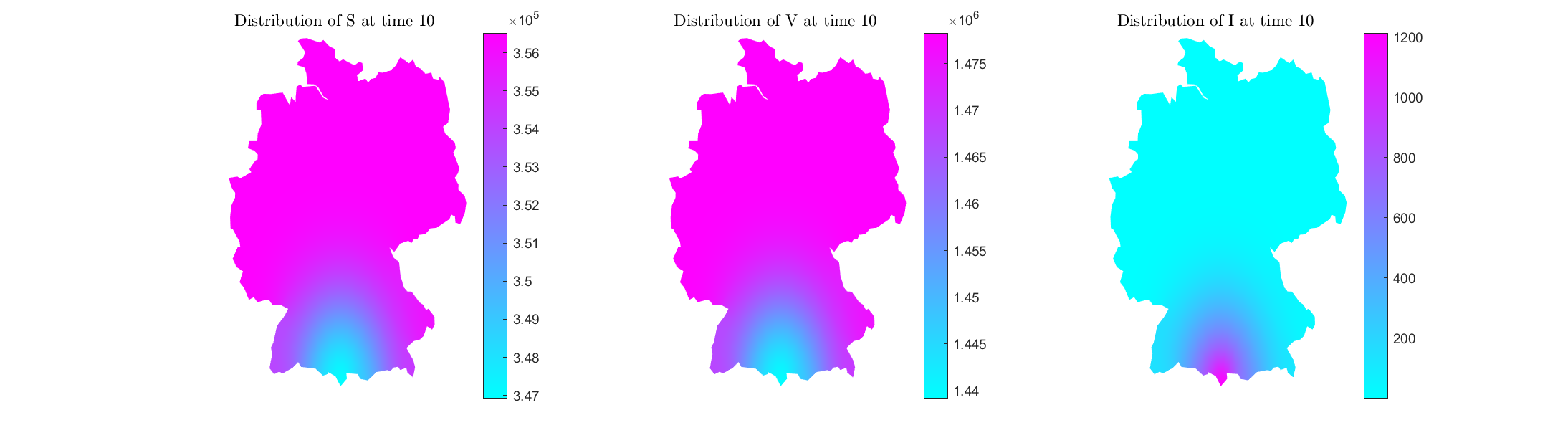}
\caption{Solution after $t_f=10$ days. \label{1peakosc10}}
\end{figure}

\begin{figure}[ht]
\includegraphics[width = \textwidth]{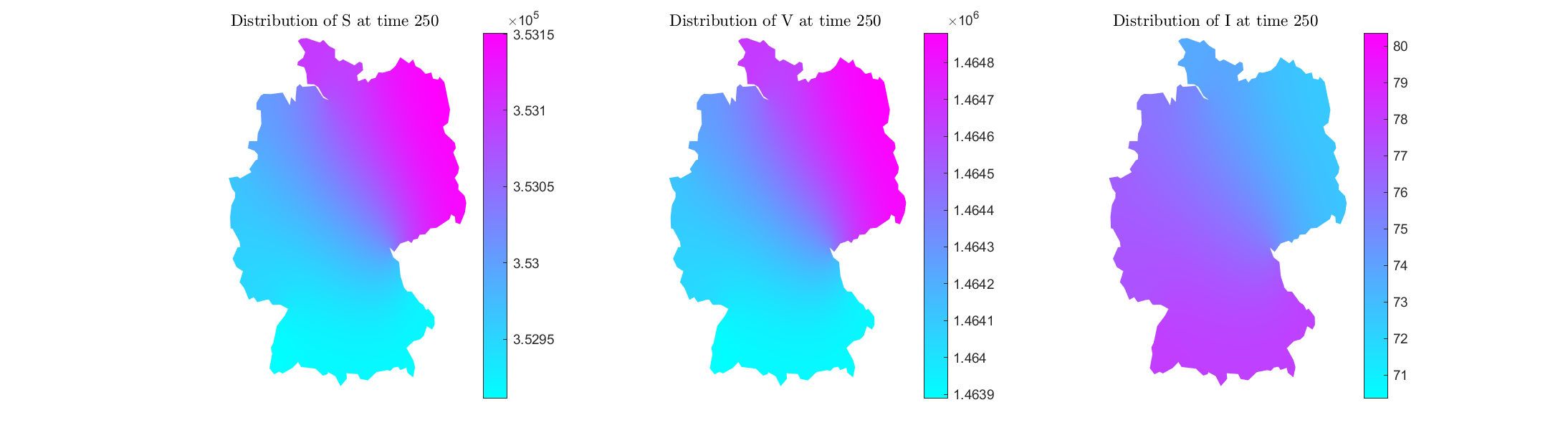}
\caption{Solution after $t_f=250$ days. \label{1peakosc250}}
\end{figure}

\begin{figure}[ht]
\includegraphics[width = \textwidth]{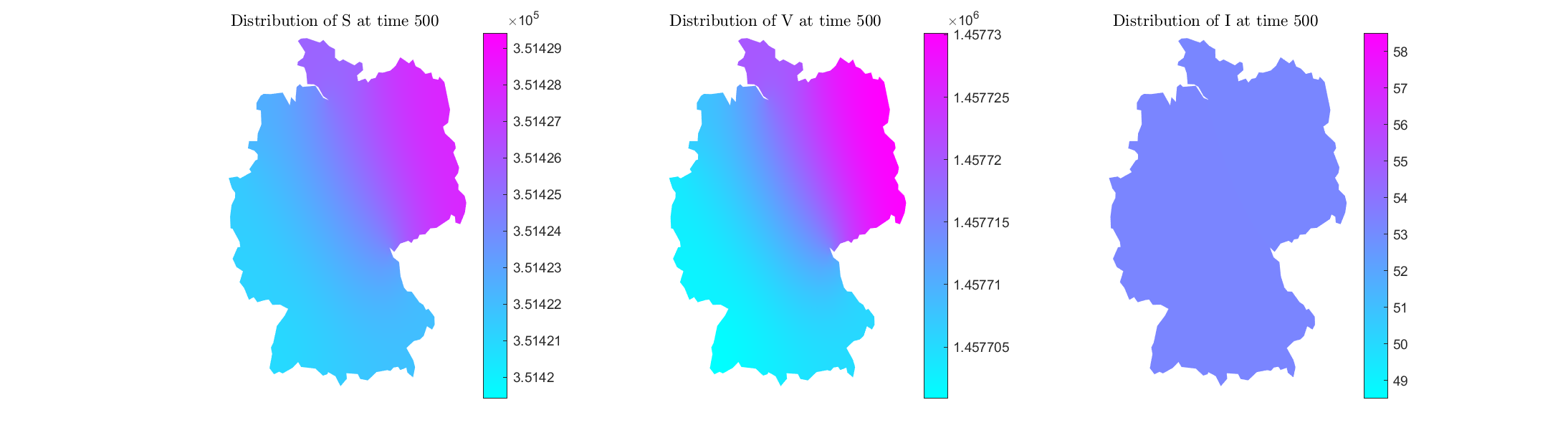}
\caption{Solution after $t_f=500$ days. \label{1peakosc500}}
\end{figure}

\begin{figure}[ht]
\includegraphics[width=0.9\textwidth]{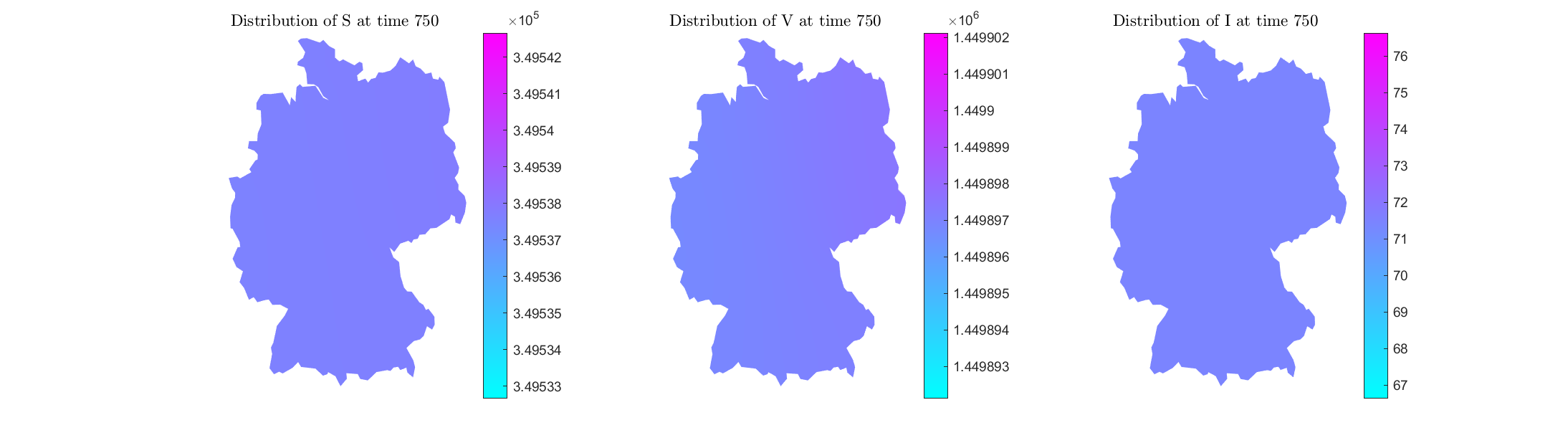}
\caption{Solution after $t_f=750$ days. The solution is almost spatially constant. \label{1peakosc750}}
\end{figure}

%\FloatBarrier
%\subsubsection{Time-dependent parameters}
%\label{Time-dependent contact rate}
Temperature has been observed to play an important role in the spread of Covid-19 in Germany \cite{ganegoda2021interrelationship}. To take this into account, we assume that the transmission rate coefficient $\beta$ is time dependent. For this case, we choose $\beta(t)=\beta_0\left(1+\alpha\cos(\frac{2\pi*t}{12}+b)\right)$, where $\beta_0=0.92429$ (see table \ref{tab_param_germany_Vacc}), $\alpha=0.0228$ and $b=\pi/6$ for which there is $\mathcal R_c =1.2058$. Figures \ref{1peakosc10}-\ref{1peakosc750} illustrate the distribution of the state of the model when the final time $t_{final}$ is equal to $10$, $250$, $500$ and $750$, respectively. The initial values are the same as in Figure \ref{1PeaksInit}. Note that for $t_{final}=750$ days (see figure \ref{1peakosc750}), the transmission rate $\beta$ is equal to $0.85445$, which corresponds to a value of the control reproduction number greater than one ($\mathcal{R}_c=1.04228$); while for $t_{final}=1000$ days, we obtain $\beta$ equal to $0.81946$ less than $\beta_{crit}$, and which corresponds to a value of the control reproduction number less than one ($\mathcal{R}_c=0.99960$).

\begin{figure}[ht]
\includegraphics[width=0.9\textwidth]{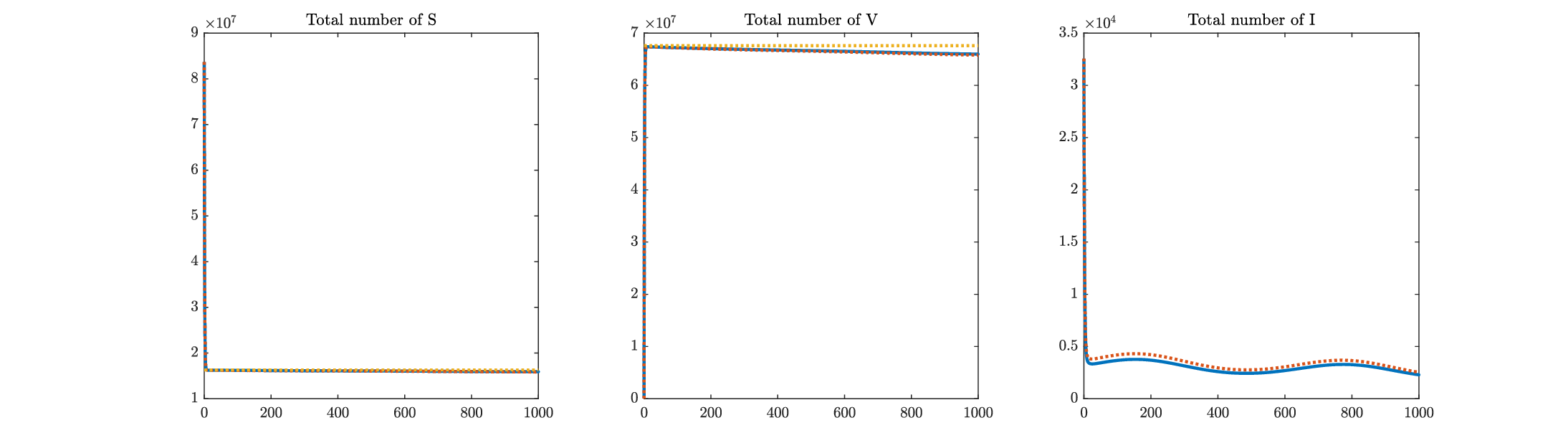}
\caption{Total number of members in each compartment over time (blue), compared to the ODE model (red) and the disease-free state (yellow). All parameter values are ones of Table \ref{tab_param_germany_Vacc} except
$\beta(t)=\beta_0(1+\alpha\cos((\dfrac{2\pi t}{12})+b)$, where, where $\beta_{crit} = 0.819792192423568$ is the critical value of $\beta$ for which there is $\mathcal R_c = 1$. \label{1peakosctotal}}
%\end{center}
\end{figure}

The total number of members in each compartment over time (blue), compared to the ODE model (red) and the disease-free state (yellow) when
$\beta(t)=\beta_0(1+0.0228\cos((\frac{2\pi t}{12})+\frac{\pi}{6})$ is depicted in figure \ref{1peakosctotal}. It gives the impression that the number of infected individuals remains bounded, in contrast to the situation in Figure \ref{1PeaksTotal} ($\beta$ constant, greater than
$\beta_{crit}$).

\FloatBarrier
\subsubsection{Cameroon}
In this part, we turn to the simulation of the Covid-19 situation in Cameroon. For our simulations, we used the initial coefficients and values of table \ref{tab_param_Cam_Vacc}. We present the results for two initial situations, one in which all infectious are concentrated in the very south, and the other in which there is an additional center of infection in the very north of the country.
	
\begin{figure}[ht]
	%\begin{center}
	%\includegraphics[scale=0.5]{InitALL}
	\includegraphics[width = \textwidth]{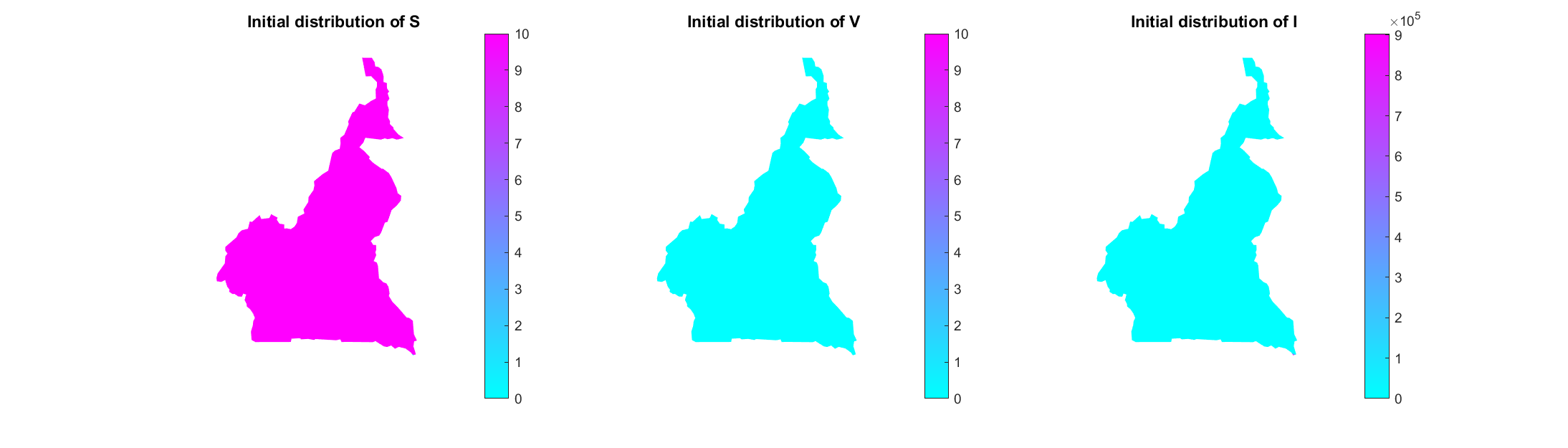}
	\caption{Initial state. \label{InitialK}}
	%\end{center}
\end{figure}

\begin{figure}[ht]
	%\begin{center}
	%\includegraphics[scale=0.5]{InitALL}
	\includegraphics[width = \textwidth]{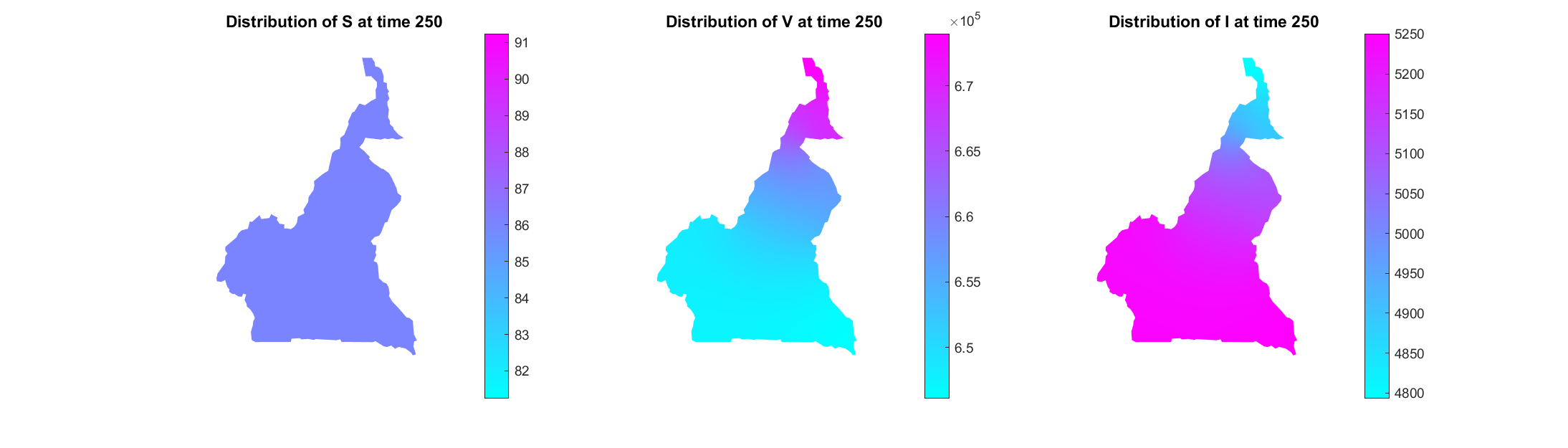}
	\caption{Solution after $t_f=250$ days.\label{1peak250K}}
	%\end{center}
\end{figure}

\begin{figure}[ht]
	%\begin{center}
	%\includegraphics[scale=0.5]{InitALL}
	\includegraphics[width = \textwidth]{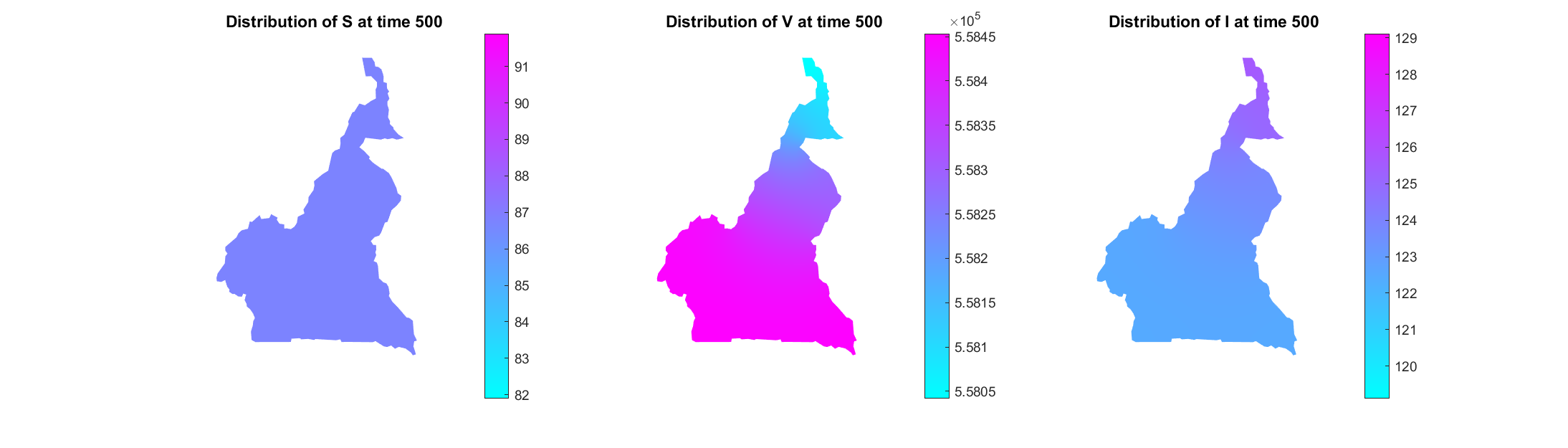}
	\caption{Solution after $t_f=500$ days.\label{1peak500K}}
	%\end{center}
\end{figure}

\begin{figure}[ht]
	%\begin{center}
	%\includegraphics[scale=0.5]{InitALL}
	\includegraphics[width = \textwidth]{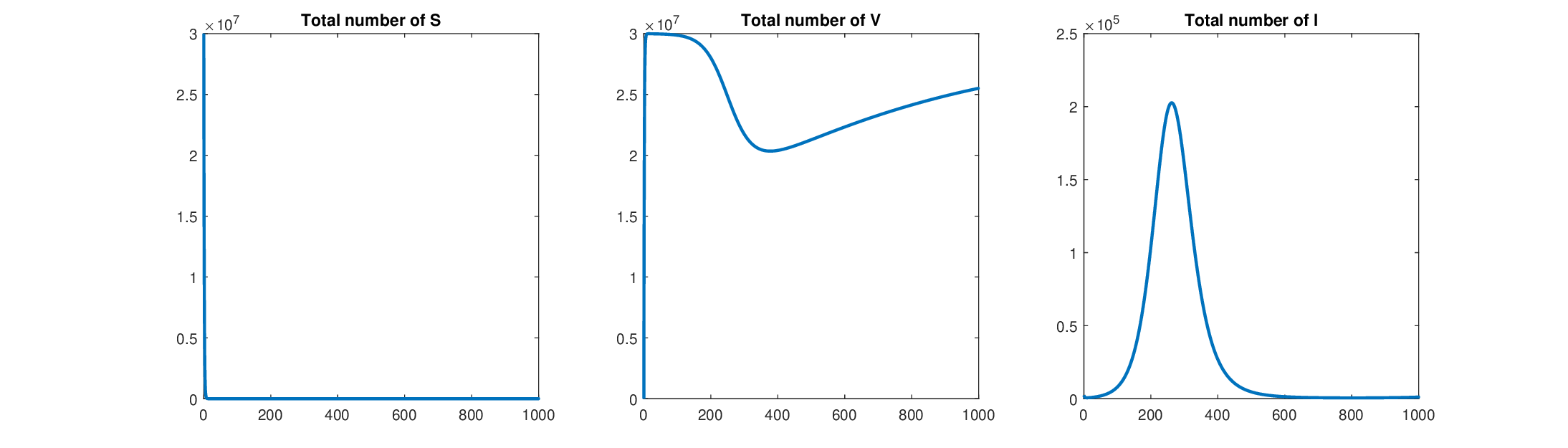}
	\caption{Final state of the solution.\label{1peakFinal_totalK}}
	%\end{center}
\end{figure}

	\paragraph{One peak}
	For this simulation, we assume that all infective individuals are concentrated in Cameroon's southern corner. As in the previous section, we present the results for $\mathfrak{S}$, $\mathfrak{V}$ and $\mathfrak{I}$ for various intermediate times and the total number of those compartments over the whole time span. The results are displayed in figures
\ref{InitialK} to \ref{1peakFinal_totalK}.
%\begin{figure}[ht]
	%\begin{center}
	%\includegraphics[width = \textwidth]{1peak_10_K}
	%\caption{Solution after $t_f=10$ days. \label{1peak10K}}
	%\end{center}
%\end{figure}
%\begin{figure}[ht]
	%\begin{center}
	%\includegraphics[scale=0.5]{InitALL}
	%\includegraphics[width = \textwidth]{1peak_750_K}
	%\caption{Solution after $t_f=750$ days.\label{1peak750K}}
	%\end{center}
%\end{figure}
%\begin{figure}[ht]
	%\begin{center}
	%\includegraphics[scale=0.5]{InitALL}
	%\includegraphics[width = \textwidth]{1peak_final_K}
	%\caption{Final state of the solution.\label{1peakFinalK}}
	%\end{center}
%\end{figure}

\begin{figure}[ht]
	%\begin{center}
	%\includegraphics[scale=0.5]{InitALL}
	\includegraphics[width = \textwidth]{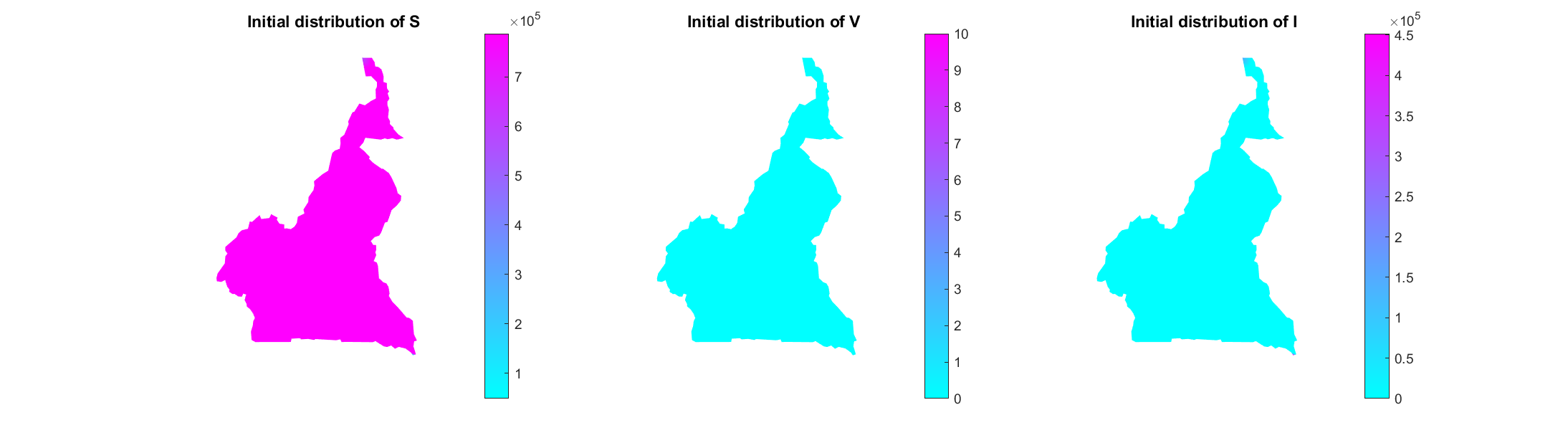}
	\caption{Initial state. \label{2peakInitK}}
	%\end{center}
\end{figure}

\begin{figure}[ht]
	%\begin{center}
	%\includegraphics[scale=0.5]{InitALL}
	\includegraphics[width = \textwidth]{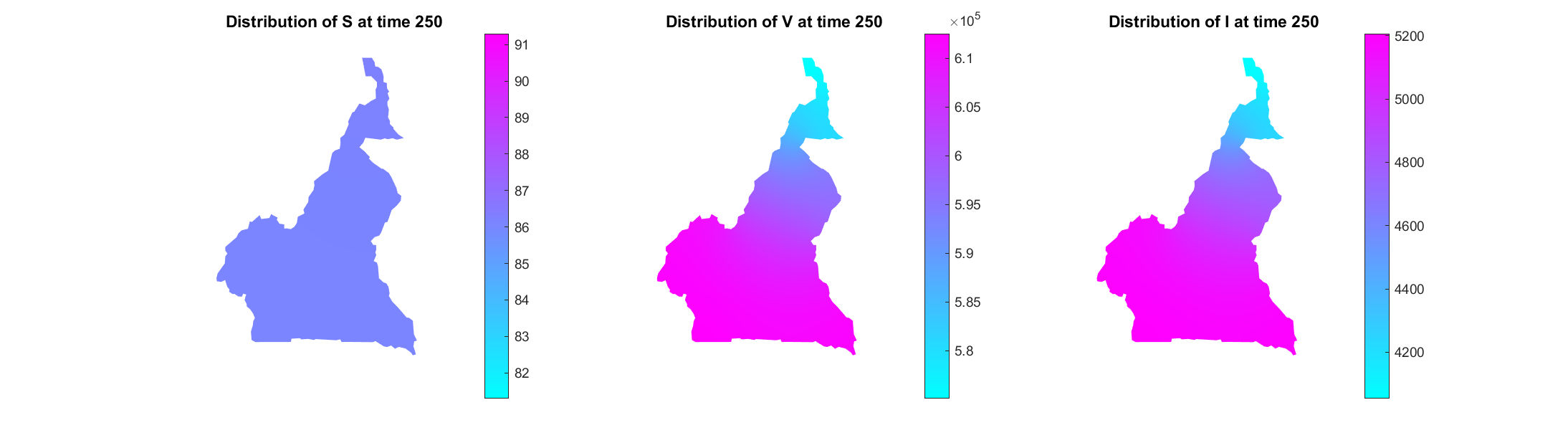}
	\caption{Solution after $t_f=250$ days.\label{2peak250K}}
	%\end{center}
\end{figure}

\begin{figure}[ht]
	%\begin{center}
	%\includegraphics[scale=0.5]{InitALL}
	\includegraphics[width = \textwidth]{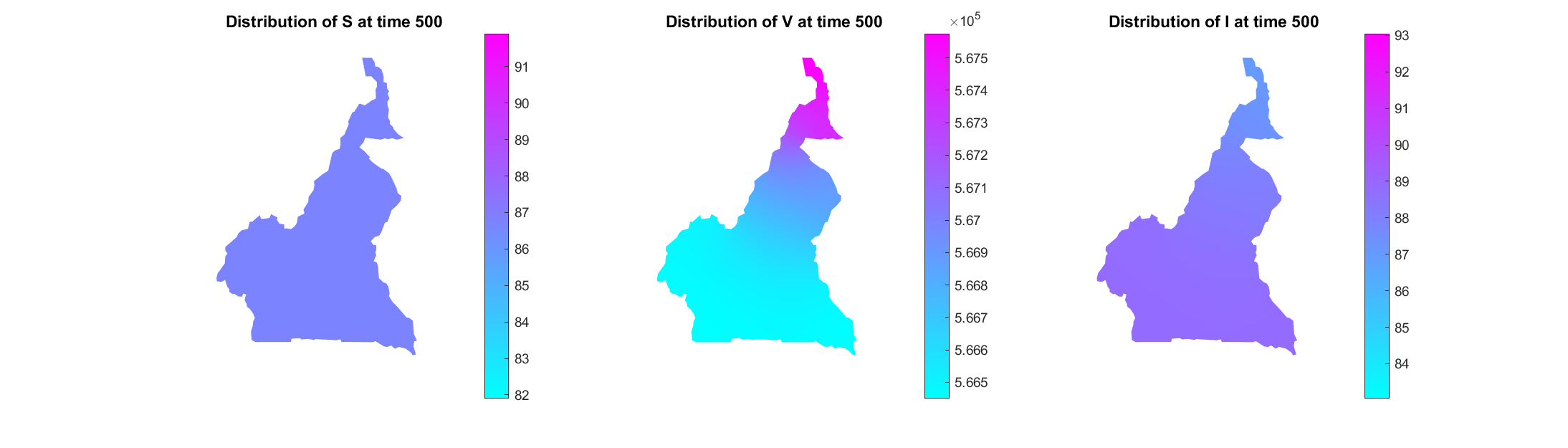}
	\caption{Solution after $t_f=500$ days.\label{2peak500K}}
	%\end{center}
\end{figure}

\begin{figure}[ht]
	%\begin{center}
	%\includegraphics[scale=0.5]{InitALL}
	\includegraphics[width = \textwidth]{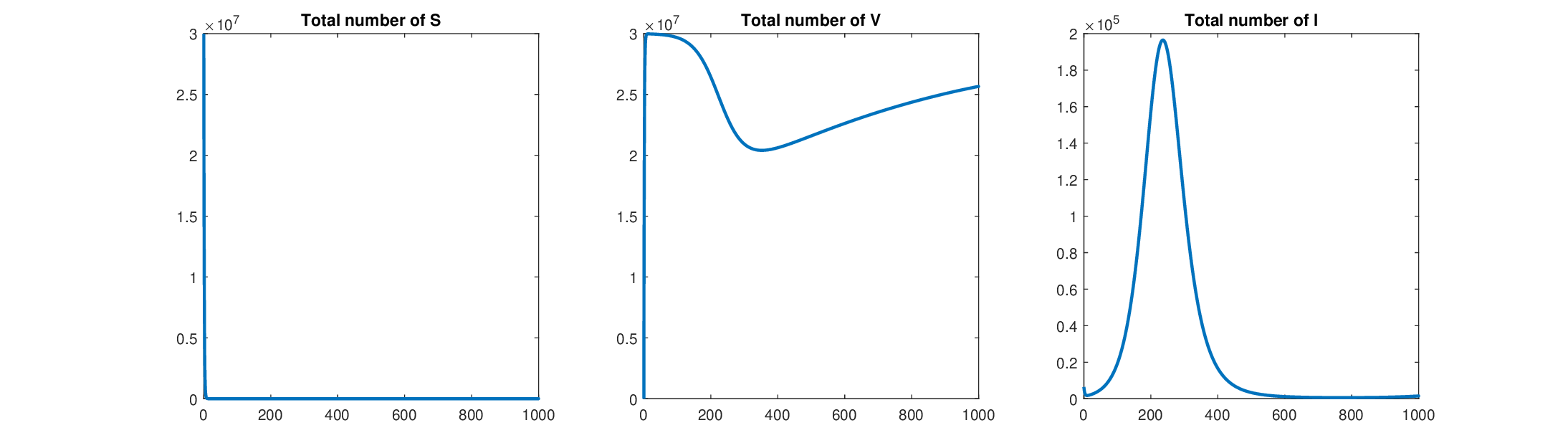}
	\caption{Final state of the solution.\label{2peakFinalTotalK}}
	%\end{center}
\end{figure}

%\FloatBarrier
\paragraph{Two peaks}
For this situation, we assume that the infective individuals are in the very north or the very south of the country. For small times, the figures for this situation are clearly distinguishable from those in the previous section, but after a while the diffusion process takes over and eliminates any spatial inhomogenity of the solution. The results are displayed in figures \ref{2peakInitK} to \ref{2peakFinalTotalK}
%\begin{figure}[ht]
	%\begin{center}
	%\includegraphics[scale=0.5]{InitALL}
	%\includegraphics[width = \textwidth]{2peak_10_K}
	%\caption{Solution after $t_f=10$ days.\label{2peak10K}}
	%\end{center}
%\end{figure}
%\begin{figure}[ht]
	%\begin{center}
	%\includegraphics[scale=0.5]{InitALL}
	%\includegraphics[width = \textwidth]{2peak_750_K}
	%\caption{Solution after $t_f=750$ days.\label{2peak750K}}
	%\end{center}
%\end{figure}
%\begin{figure}[ht]
	%\begin{center}
	%\includegraphics[scale=0.5]{InitALL}
	%\includegraphics[width = \textwidth]{2peak_final_K}
	%\caption{Final state of the solution.\label{2peakFinalK}}
	%\end{center}
%\end{figure}

\begin{figure}[ht]
	%\begin{center}
	\includegraphics[width = \textwidth]{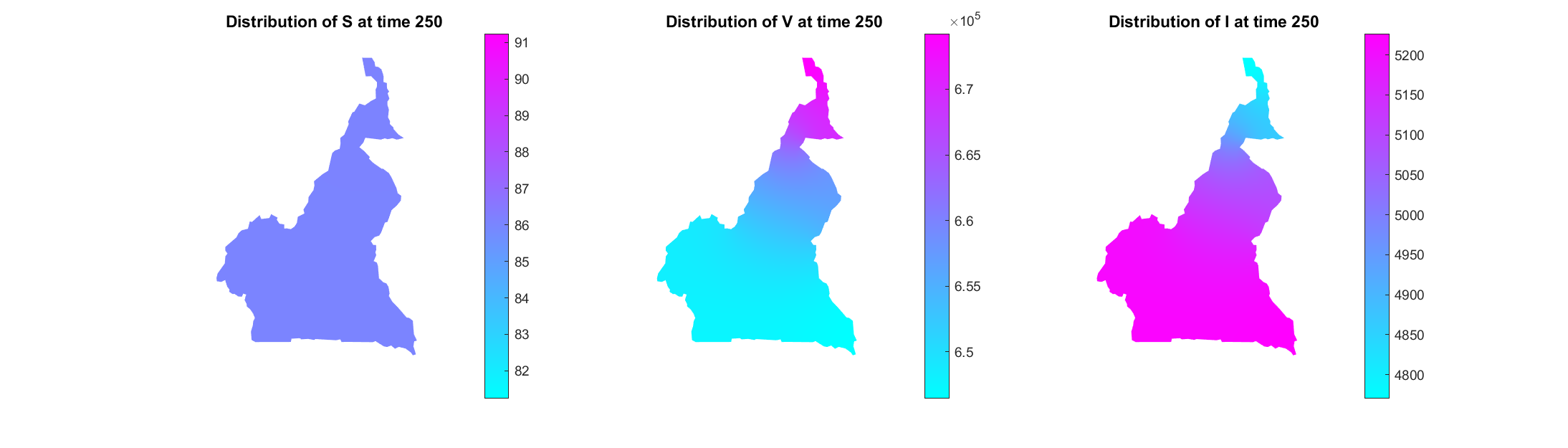}
	\caption{State of the solution with oscillations for $t_f=250$ days.\label{oscill250K}}
	%\end{center}
\end{figure}

\begin{figure}[ht]
	%\begin{center}
	\includegraphics[width = \textwidth]{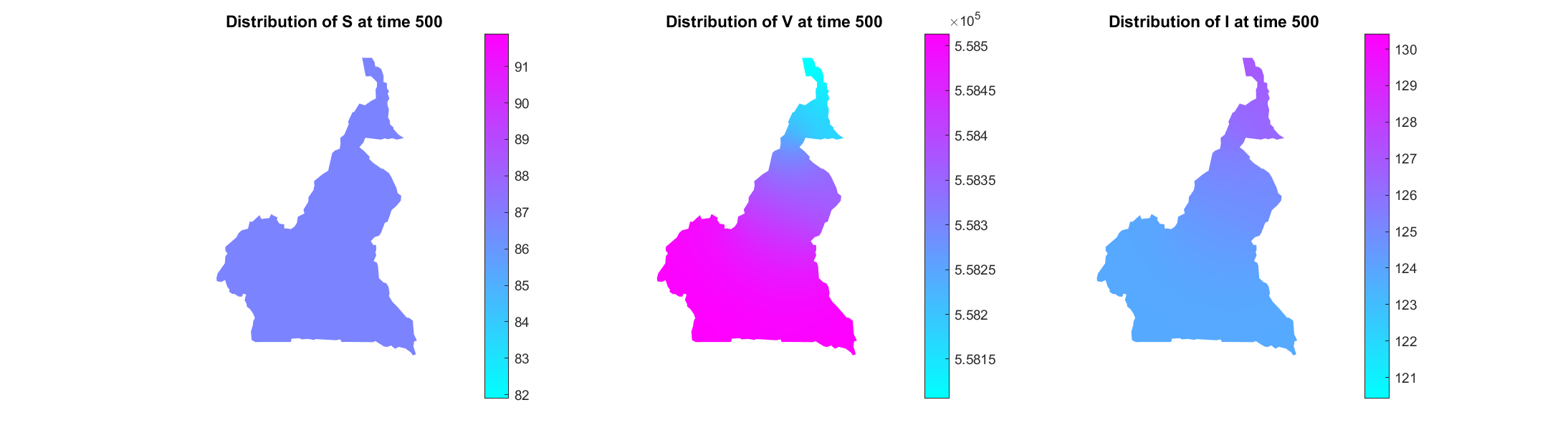}
	\caption{State of the solution with oscillations for $t_f=500$ days.\label{oscill500K}}
	%\end{center}
\end{figure}

\paragraph{Time-dependent contact rate}
In what follows, we suppose that
$\beta(t)=\beta_0(1+\alpha\cos((\dfrac{2\pi t}{12})+b)$, where
$\beta_0=0.399092568990682$ (see table \ref{tab_param_Cam_Vacc}), $\alpha=0.0228$ and $b=\pi/6$. The results are displayed in figures \ref{oscill250K}-\ref{oscill500K}.
%\begin{figure}[ht]
	%\begin{center}
	%\includegraphics[width = \textwidth]{1peak_osc_10_K}
	%\caption{Final state of the solution with oscillations for $t_f=10$ days.\label{oscill10K}}
	%\end{center}
%\end{figure}

%\begin{figure}[ht]
	%\begin{center}
	%\includegraphics[width = \textwidth]{1peak_osc_750_K}
	%\caption{State of the solution with oscillations for $t_f=750$ days. \label{oscill750K}}
	%\end{center}
%\end{figure}

%\newpage

\FloatBarrier
\section{Conclusion and perspectives}
In this work, we formulated two Covid-19 transmission dynamics models, first using ordinary differential equations, and secondly, partial differential equations (reaction-diffusion model). For both models, we computed the control reproduction number and proved the global asymptotic stability of the disease-free equilibrium point whenever the control reproduction number is less than one. We also proved that both models admit at least one endemic equilibrium point whenever the control reproduction number is greater than one, which implies that the disease-free equilibrium becomes unstable. Using the daily infected reported cases in Germany (resp. in Cameroon) from December 31, 2020, to February 28, 2021 (resp from April 2021 to April 2022), we calibrated the ODE model by estimated model parameters. We found that the control reproduction number, $\mathcal{R}_c$, is approximately equal to $1.13$ for Germany (resp. $1.2554$ for Cameroon) which confirms that, even if the vaccination level is high, Covid-19 will be present in these countries, and this for the next years.

The final part of the work concerns numerical studies. After presenting the numerical method used to simulate the models, we performed various numerical simulations to validate our theoretical results. Indeed, several cases are considered: First, of all, we considered the case where all parameters are constant, with different final times , followed by the case where the transmission rate coefficient $\beta$ is time dependent. Third, we considered the case in which the initial population is entirely susceptible to infection, except for one small region in the very south of Germany, where there are also infected and exposed \ldots people, and the case in which we added a second peak in western Germany. These two cases were followed by the case in which only the state of Bavaria holds exposed, asymptomatic infected, and symptomatic infected individuals. For each of the above cases, we compared the ODE model with the PDE model by drawing in the same panel the total numbers in each compartment for the PDE and ODE model, as well as the corresponding disease-free equilibrium points (yellow). In order to directly compare the two models, we chose spatially constant parameters and initial values for the PDE to rule out diffusion effects. This permits us to conclude that the total numbers in each compartment for the PDE and ODE models coincide in this case, and from a quantitative point of view, the ODE model and the PDE model then give the same results.

Although these two countries (Germany and Cameroon) have practically the same number of base reproductions, the damage caused by COVID-19 in these countries is different in terms of loss of human life. In fact, the predictions of an epidemic in the sub-Saharian Africa countries by the World Health Organization have not been confirmed\footnote{\url{https://www.afro.who.int/fr/node/12206}}. Unlike the countries of Western Europe, which have the financial resources and health centers with better technical facilities and an excellent qualified workforce, the outbreak of the disease in Cameroon has been less severe. Some virus specialists have claimed that taking chloroquine to treat malaria has played an important role in immunizing the population against COVID-19. It is also possible that the use of drugs derived from traditional research has boosted the population's immune system\footnote{\url{https://international.la-croix.com/fr/afrique/covid-19-le-remede-de-mgr-kleda-autorise-par-le-gouvernement-camerounais}}.

In the present study, we did not take into account the fact that model parameters can depend on time and space. Indeed, the transmission rate $\beta$, for example, should not be the same in a country as Germany which has sixteen federal states with different population sizes and densities. Thus, estimating some model parameters for each German state, taking into account population movement between each state, constitutes a direct perspective of this work.

\section*{Acknowledgments}
The authors thank the Ministry of Science, Research and the Arts of the State of Baden-W\"urttemberg (Ministerium f\"ur Wissenschaft, Forschung und Kunst Baden-W\"urttemberg) for the grant supporting a joint research project within the 2022 initiative "Science Cooperation Africa (2021/2022)" which allowed the first author to have a research stay in the Department of Mathematics and Statistics of the University of Konstanz, Germany. The first author also thanks the Zukunftskolleg for providing him with all the necessary logistics during this research stay at the University of Konstanz.

\section*{Conflict of interest}
The authors declare that they have no conflict of interest.

%\section*{Declarations}
%\begin{itemize}
%\item The authors declare that they have no conflict of interest.
%\item Not applicable
%\item Consent for publication
%\item Not applicable
%\item Not applicable
%\item Not applicable
%\item All authors contribute equally to to this work
%\end{itemize}

%\bibliographystyle{spmpsci}
\bibliographystyle{unsrt}
% inclusion de la biblio
\bibliography{Bibliographie}
%\end{linenumbers}
\end{document}